\documentclass[12pt,1p] {elsarticle}
\usepackage{amsmath}
\usepackage{amsfonts}
\usepackage{amssymb}
\usepackage{color}
\usepackage{graphicx}
\usepackage{local}
\usepackage{theorem}
\usepackage{alg}
\usepackage{subfigure}
\usepackage[usenames,dvipsnames]{xcolor}
\usepackage{tikz}
\usepackage{pgfplots}
\usetikzlibrary{shapes}
\usepackage{multirow}
\tikzset{
  amrpatch/.is family,
  amrpatch,
  shiftx/.initial=0,
  shifty/.initial=0,
  shiftz/.initial=0,
  dimx/.initial=3,
  dimy/.initial=3,
  dimz/.initial=3,
  scale/.initial=1,
  densityx/.initial=1,
  densityy/.initial=1,
  densityz/.initial=1,
  rotation/.initial=0,
  anglex/.initial=0,
  angley/.initial=90,
  anglez/.initial=210,
  scalex/.initial=1,
  scaley/.initial=1,
  scalez/.initial=0.8,
  front/.style={draw=black,fill=white},
  top/.style={draw=black,fill=white},
  right/.style={draw=black,fill=white},
  edgestyle/.style ={line width=0pt},
  backedgestyle/.style ={line width=1pt},
  frontedgestyle/.style={line width=1pt},
  innerlinestyle/.style={line width=1pt},
  fillopacity/.style={fill opacity=1},
}

\newcommand{\tikzamrpatchkey}[1]{\pgfkeysvalueof{/tikz/amrpatch/#1}}

\newcommand{\amrpatch}[1]{
   \tikzset{amrpatch,#1} % Process Keys passed to command
   \pgfmathsetlengthmacro{\vectorxx}{\tikzamrpatchkey{scalex}*cos(\tikzamrpatchkey{anglex})*28.452756}
   \pgfmathsetlengthmacro{\vectorxy}{\tikzamrpatchkey{scalex}*sin(\tikzamrpatchkey{anglex})*28.452756}
   \pgfmathsetlengthmacro{\vectoryx}{\tikzamrpatchkey{scaley}*cos(\tikzamrpatchkey{angley})*28.452756}
   \pgfmathsetlengthmacro{\vectoryy}{\tikzamrpatchkey{scaley}*sin(\tikzamrpatchkey{angley})*28.452756}
   \pgfmathsetlengthmacro{\vectorzx}{\tikzamrpatchkey{scalez}*cos(\tikzamrpatchkey{anglez})*28.452756}
   \pgfmathsetlengthmacro{\vectorzy}{\tikzamrpatchkey{scalez}*sin(\tikzamrpatchkey{anglez})*28.452756}
   \begin{scope}[scale=\tikzamrpatchkey{scale}, rotate=\tikzamrpatchkey{rotation},
                 x={(\vectorxx,\vectorxy)}, y={(\vectoryx,\vectoryy)}, z={(\vectorzx,\vectorzy)}]
   \newcommand{\dimx}{\tikzamrpatchkey{dimx}}
   \newcommand{\dimy}{\tikzamrpatchkey{dimy}}
   \newcommand{\dimz}{\tikzamrpatchkey{dimz}}
   \newcommand{\shiftx}{\tikzamrpatchkey{shiftx}}
   \newcommand{\shifty}{\tikzamrpatchkey{shifty}}
   \newcommand{\shiftz}{\tikzamrpatchkey{shiftz}}
   \pgfmathsetmacro{\steppingx}{1/\tikzamrpatchkey{densityx}}
   \pgfmathsetmacro{\steppingy}{1/\tikzamrpatchkey{densityy}}
   \pgfmathsetmacro{\steppingz}{1/\tikzamrpatchkey{densityz}}
   \pgfmathsetmacro{\secondx}{2*\steppingx}
   \pgfmathsetmacro{\secondy}{2*\steppingy}
   \pgfmathsetmacro{\secondz}{2*\steppingz}

   \ifthenelse{\lengthtest{\secondx pt > \dimx pt}}
   {
      \pgfmathsetmacro{\secondx}{\dimx+0.00001}
   }
   {
      \pgfmathsetmacro{\secondx}{2*\steppingx}
   }
   \ifthenelse{\lengthtest{\secondy pt > \dimy pt}}
   {
      \pgfmathsetmacro{\secondy}{\dimy+0.00001}
   }
   {
      \pgfmathsetmacro{\secondy}{2*\steppingy}
   }
   \ifthenelse{\lengthtest{\secondz pt > \dimz pt}}
   {
      \pgfmathsetmacro{\secondz}{\dimz+0.00001}
   }
   {
      \pgfmathsetmacro{\secondz}{2*\steppingz}
   }

   \foreach \x in {\steppingx,\secondx,...,\dimx}
   {
      \foreach \y in {\steppingy,\secondy,...,\dimy}
      {
         \pgfmathsetmacro{\lowx}{(\x-\steppingx)}
         \pgfmathsetmacro{\lowy}{(\y-\steppingy)}
         \filldraw[amrpatch/front,amrpatch/edgestyle,amrpatch/fillopacity] (\lowx+\shiftx,\lowy+\shifty,\dimz+\shiftz) -- (\lowx+\shiftx,\y+\shifty,\dimz+\shiftz) -- (\x+\shiftx,\y+\shifty,\dimz+\shiftz) -- (\x+\shiftx,\lowy+\shifty,\dimz+\shiftz) -- cycle;
      }
   }

   \foreach \x in {\steppingx,\secondx,...,\dimx}
   {
      \foreach \z in {\steppingz,\secondz,...,\dimz}
      {
         \pgfmathsetmacro{\lowx}{(\x-\steppingx)}
         \pgfmathsetmacro{\lowz}{(\z-\steppingz)}
         \filldraw[amrpatch/top,amrpatch/edgestyle,amrpatch/fillopacity] (\lowx+\shiftx,\dimy+\shifty,\lowz+\shiftz) -- (\lowx+\shiftx,\dimy+\shifty,\z+\shiftz) -- (\x+\shiftx,\dimy+\shifty,\z+\shiftz) -- (\x+\shiftx,\dimy+\shifty,\lowz+\shiftz) -- cycle;
      }
   }
   \foreach \y in {\steppingy,\secondy,...,\dimy}
   {
      \foreach \z in {\steppingz,\secondz,...,\dimz}
      {
         \pgfmathsetmacro{\lowy}{(\y-\steppingy)}
         \pgfmathsetmacro{\lowz}{(\z-\steppingz)}
         \filldraw[amrpatch/right,amrpatch/edgestyle,amrpatch/fillopacity] (\dimx+\shiftx,\lowy+\shifty,\lowz+\shiftz) -- (\dimx+\shiftx,\lowy+\shifty,\z+\shiftz) -- (\dimx+\shiftx,\y+\shifty,\z+\shiftz) -- (\dimx+\shiftx,\y+\shifty,\lowz+\shiftz) -- cycle;
      }
   }

   \end{scope}
}

\pgfplotscreateplotcyclelist{bigmark}
{
   {black,mark=triangle,mark size=4pt},
   {black,mark=|,       mark size=6pt},
   {black,mark=square,  mark size=4pt},
   {black,mark=o,       mark size=4pt},
}

\makeatother

\journal{Journal of Computational Physics}

\begin{document}

\begin{frontmatter}

\title{Dynamic Implicit 3D Adaptive Mesh Refinement for Non-Equilibrium Radiation Diffusion\tnoteref{tnote}}

\tnotetext[tnote]{Notice: This manuscript has been authored by UT-Battelle, LLC, under Contract No. DE-AC05-00OR22725 with
 the U.S. Department of Energy. The United States Government retains and the publisher, by accepting the
 article for publication, acknowledges that the United States Government retains a non-exclusive,
 paid-up, irrevocable, world-wide license to publish or reproduce the published form of this
 manuscript, or allow others to do so, for United States Government purposes.
}

\author[csm]{B. Philip\corref{cor1}}
%\ead{bphilip.kondekeril@gmail.com}
\ead{bphilip.kondekeril@gmail.com}

\author[nccs]{Z. Wang}
\ead{wangz@ornl.gov}

\author[csm]{M. A. Berrill}
\ead{berrillma@ornl.gov}

\author[nccs]{M. Rodriguez Rodriguez}
\ead{manuro@live.de}

\author[inl]{M. Pernice}
\ead{michael.pernice@inl.gov}

\address[csm]{
Computer Science and Mathematics Division\\
Oak Ridge National Laboratory\\
Oak Ridge TN 37831-6164}

\address[nccs]{
National Center for Computational Sciences\\
Oak Ridge National Laboratory\\
Oak Ridge TN 37831-6164}

\address[inl]{
Applied Computing and Visualization\\
Energy, Environment Science and Technology \\
Idaho National Laboratory\\
Idaho Falls, ID 83415-3550\\
}

\cortext[cor1]{Corresponding author}

\begin{abstract}
The time dependent non-equilibrium radiation diffusion equations are important for solving the transport of energy through radiation in optically thick regimes and find applications in several fields including astrophysics and inertial confinement fusion. The associated initial boundary value problems that are encountered often exhibit a wide range of scales in space and time and are extremely challenging to solve. To efficiently and accurately simulate these systems we describe our research on combining techniques that will also find use more broadly for long term time integration of nonlinear multiphysics systems: implicit time integration for efficient long term time integration of stiff multiphysics systems, local control theory based step size control to minimize the required global number of time steps while controlling accuracy, dynamic 3D adaptive mesh refinement (AMR) to minimize memory and computational costs, Jacobian Free Newton-Krylov methods on AMR grids for efficient nonlinear solution, and optimal multilevel preconditioner components that provide level independent solver convergence. 
\end{abstract}

\end{frontmatter}

\section{Introduction}
In the fields of astrophysics and inertial confinement fusion the time dependent non-equilibrium radiation diffusion equations are important for solving the transport of energy through radiation in an optically thick regime.  In this paper we employ a form of the model that has a flux-limited diffusion approximation (gray approximation) for the energy density coupled to a material temperature equation that incorporates a nonlinear material conduction term \cite{MihalasWeibel-Mihalas99,Olson-jcp-07,mavriplis01,Mousseau00,mousseau-jcp-03}.  This nonlinear, coupled, time dependent set of partial differential equations (PDEs) exhibits multiple temporal and spatial scales, and the associated initial boundary value problems are highly stiff and challenging to solve. As a result, they are also an excellent testbed for the development of simulation methods for long term time integration of stiff multi-physics systems.

In this paper we will limit our scope to fully implicit time integration methods. This then enables the use of timestep control methods based on accuracy considerations and enables us to leverage the theoretical advances for accuracy based timestep control that exist in the field of ordinary differential equations (ODEs). We experiment with different adaptive timestep control methods including control theory based approaches that attempt to monitor and control the temporal accuracy at each timestep and minimize the total number of timesteps required over the course of the simulation. The use of control theoretic approaches  to timestep control is new for radiation-diffusion calculations and is only beginning to be used for multi-physics calculations. Variable step fully implicit time integration is combined with 3D dynamic structured adaptive mesh refinement (AMR)\cite{BergerThesisAMR} with an objective towards minimizing the total number of degrees of freedom required over the course of a simulation. Care is however required in combining these techniques as spatial regridding during dynamic AMR can introduce non-stiff transient errors that significantly affect the behavior of timestep control algorithms and can lead to a dramatic increase in the total number of timesteps required over uniform spatial grid calculations when not properly controlled. We will report on our experiences with different timestep controllers and the modifications required in this context for AMR as the literature on this topic, particularly for multi-physics simulations, is sparse. Fully implicit time integration methods require highly efficient solution of the nonlinear systems at each timestep in order to be competitive with other methods. Here, we choose to use Jacobian Free Newton-Krylov (JFNK) methods with physics based preconditioning. JFNK methods allow us to avoid the formation of the full Jacobian matrices across AMR grid hierarchies which can be problematic and programming intensive for flux based finite volume  discretizations on 3D AMR grid hierarchies which incorporate coarse-fine interpolation across grid levels.  JFNK methods often obtain their efficiency from careful design of preconditioners. Efficient preconditioners on uniform grids for JFNK methods often employ multigrid solvers to tackle elliptic components to deliver grid independent performance. In the context of AMR, particularly for problems with elliptic components, preconditioner performance can degrade as the number of refinement levels in the AMR hierarchy increases if proper care is not paid to coupling between levels.  By employing suitable multilevel preconditioner components we will demonstrate level independent performance of our nonlinear solvers for non-equilibrium radiation diffusion applications.

The remainder of this paper is organized as follows. Section 2 of this paper surveys related work in the context of equilibrium and non-equilibrium radiation diffusion problems. Section 3 describes the model problem and its temporal and spatial discretization. Section 4 describes the JFNK method and the multilevel preconditioners employed. Section 5 presents numerical results and Section 6 presents conclusions and directions for future and ongoing work.

\section{Related work}
In \cite{knoll-jqsrt-99}, Rider, Knoll and Olson introduced the idea of physics based preconditioning in 1D for non-equilibrium radiation diffusion problems. Further work by Mousseau, Knoll, Rider \cite{Mousseau00} and Mousseau, Knoll \cite{mousseau-jcp-03} extended this methodology to problems in 2D on uniform grids. Their work related to physics-based preconditioning will be leveraged here with major extensions for 3D AMR grids and multilevel preconditioners. In \cite{mavriplis01}, Mavriplis compared two different approaches to solving the nonlinear systems at each timestep by considering Newton-Multigrid and Full Approximation Scheme (FAS) using agglomeration ideas on unstructured grids for this problem.  In \cite{Olson-jcp-07} Olson considers the use of efficient operator split time integration schemes on uniform grids. Work by Lowrie et. al.\cite{lowrie-jcp-04} compares different time integration methods for non-equilibrium radiation diffusion while Brown, Shumaker, Woodward \cite{brown-jcp-04} focus on fully implicit methods and high order time integration on uniform grids. The motivation to consider automatic timestep control in our work was partially derived from \cite{brown-jcp-04}. We build on their work to further consider the use of the control theory based timestep controllers that provide computational stability as described in \cite{Soderlind2002} and related references and consider modifications that are required for AMR. Glowinski, Toivanen \cite{glowinski-jcp-05} consider using automatic differentiation and system multigrid. Shestakov, Greenough, and Howell \cite{shestakov-jqsrt-05} consider pseudo-transient continuation on AMR grids using an alternative formulation. Also worth mentioning is related work for  {\it equilibrium} radiation diffusion problems. Stals \cite{Stals03} compares the performance of Newton-Multigrid and FAS with local refinement on \textsl{unstructured} grids and Pernice, Philip \cite{PernicePhilip06} use JFNK with a Fast Adaptive Composite Grid (FAC) preconditioner on AMR grids for single physics equilibrium radiation-diffusion on structured adaptive mesh refinement (SAMR) grids. 
\section{Problem formulation and discretization}
\subsection{Model problem}

The non-dimensional model equations considered in this paper are given by \cite{MihalasWeibel-Mihalas99,Olson-jcp-07,mavriplis01,Mousseau00,mousseau-jcp-03}:
\begin{eqnarray}
   \frac{\partial E}{\partial t} -\nabla \cdot ( D_E \nabla E)& =&\;\;\;\sigma_a(T^4-E) \qquad \mbox{ in } \Omega, \label{eqn:E_td}\\
   \frac{\partial T}{\partial t} -\nabla \cdot (D_T \nabla T) &=&-\sigma_a(T^4-E) \qquad \mbox{ in } \Omega, \label{eqn:T_td}
%\label{pdes}
\end{eqnarray}
where $E$ is the radiation energy density, $T$ the material temperature,
$ \nabla $ the gradient, $ \nabla \cdot $ the divergence operator, and
$D_E$ and $D_T$ are nonlinear diffusion coefficients given by
\begin{center}
\begin{eqnarray*}
D_E &=& \frac {1}{3\sigma_a},\\ 
D_T &=& kT^{\frac{5}{2}},
\end{eqnarray*}
\end{center}
where $\sigma_a$ is the photon absorption cross section.  $\sigma_a$ is modeled by a constitutive law of the form
\begin{equation}
\sigma_a = z^3T^{-3} \label{eqn:sigma}
\end{equation}
with $z$ being the material atomic number and $k=0.01$ in our experiments. 
$D_E$ is usually flux limited to prevent non-physical effects and we use the Wilson limiter\cite{BowersWilson91}:
\begin{center}
\begin{eqnarray*}
D_E &=& \frac {1}{\left(3\sigma_a+ \frac {\|\nabla E\|}{E}\right)}.\\ 
\end{eqnarray*}
\end{center}
For the purposes of this paper this model will suffice, though more sophisticated models suited for production level codes are considered in \cite{OlsonMorel99}. 
%Quoting from Olson and Morel \cite{OlsonMorel99}, ``.. flux-limited diffusion is simply an interpolation formula. Its accuracy anywhere other than in asymptotic limits is unknown and very problem dependent".

An initial boundary value problem (IBVP) for (\ref{eqn:E_td})-(\ref{eqn:T_td}) is posed on the unit cube domain, $\Omega = [0,1]^3$, with initial conditions
\begin{equation}
  E = E_0,  \;\;\; T=(E_{0})^{\frac{1}{4}}                   \qquad \mbox{ at } t = 0,
\label{ics}
\end{equation}
and boundary conditions
\renewcommand{\arraystretch}{1.2}
\begin{eqnarray}
  \frac{1}{2}\mathsf{n} \cdot D_r \nabla E + \frac{E}{4} = R & \mbox{ on } \partial\Omega_{\mathcal{R}}, \; t \geq 0, \label{eqn:robinbc}\\
  \mathsf{n} \cdot D_r \nabla E = 0       & \mbox{ on } \partial\Omega_{\mathcal{N}}, \; t \geq 0, \\
  \mathsf{n} \cdot \nabla T = 0     & \mbox{on } \partial\Omega, \; t \geq 0,
  \label{eqn:bcs}
\end{eqnarray}
with $\partial{\Omega} = \partial\Omega_{\mathcal{R}} \cup \partial\Omega_{\mathcal{N}} $, $ \partial\Omega_{\mathcal{R}} \cap \partial\Omega_{\mathcal{N}}=\phi$.
The Robin boundary conditions for the energy density are defined on $ \partial\Omega_{\mathcal{R}} $,
which consists of the left and right faces at $ x = 0 $ and $ x = 1 $,
while Neumann boundary conditions are imposed on all other faces.
Following \cite{knoll-jqsrt-99}, we do not  flux limit the boundary conditions in (\ref{eqn:robinbc}).

\subsection{Spatial discretization}
Let $\Omega = [x_{\mathit{l}},x_{\mathit{h}}] \times
[y_{\mathit{l}},y_{\mathit{h}}] \times [z_{\mathit{l}},z_{\mathit{h}}]$ be a rectangular computational
domain.  We discretize the domain in the $x$ direction by subdividing
$[x_{\mathit{l}},x_{\mathit{h}}]$ into $n_x$ subintervals with
centers $x_i = x_{\mathit{l}} + (i+\half)\dx$ with
$\dx=(x_{\mathit{h}}-x_{\mathit{l}})/n_x$ for $i=0,\ldots,n_x-1$.
Each subinterval has faces located at $x_{i-\half} = x_i - \dx/2$ and
$x_{i+\half} = x_i + \dx/2$.  The $y$ and $z$ directions are similarly discretized by 
subdividing $[y_{\mathit{l}},y_{\mathit{h}}]$ and $ [z_{\mathit{l}},z_{\mathit{h}}]$ into
$n_y$ and $n_z$ subintervals with grid spacings $\dy=(y_{\mathit{h}}-y_{\mathit{l}})/n_y$ and
$\dz=(z_{\mathit{h}}-z_{\mathit{l}})/n_z$ respectively.  The tensor product of these subintervals partitions $\Omega$ into a
collection of computational cells $\Omega^h = \{\Omega_{i,j,k}\}$ each
of size $\dx \times \dy \times \dz$ centered at coordinates $(x_i,y_j, z_k)$.  These
ideas are readily extended to the case where $\Omega$ is a union of
non-overlapping rectangular regions, and we continue to use the same
notation $\Omega^h$ to denote such a collection of computational
cells.  Such \textsl{regular grids} are in widespread use in
computational science and engineering, and high
quality software that is tuned to regular grids, such as software for geometric
multigrid methods, is available.

\begin{figure}[ht]
\begin{center}
\includegraphics{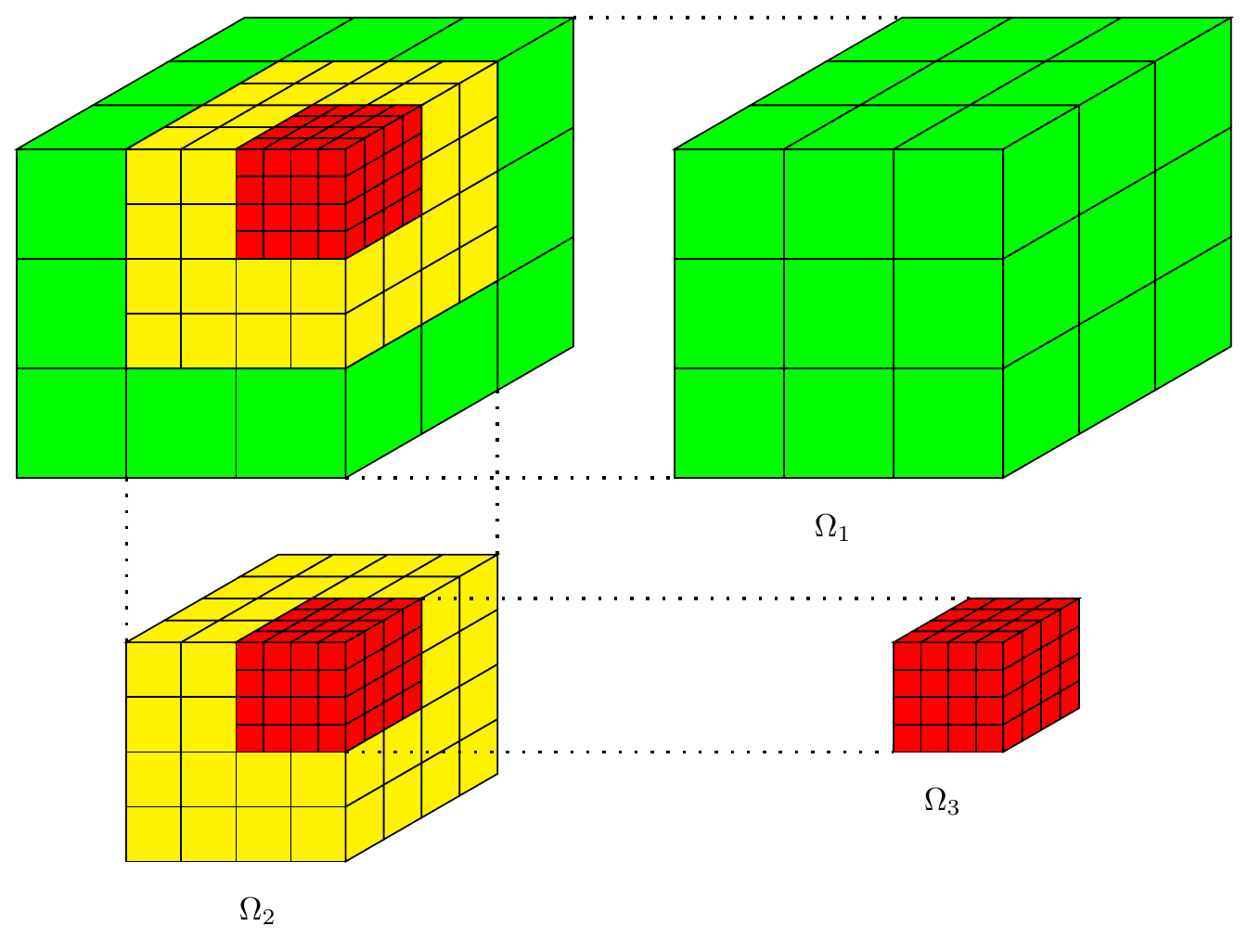}
\caption{Example of a multilevel SAMR grid with three levels.}
\label{fig:hierarchy}
\end{center}
\end{figure}

We now describe how the above approach can be extended to construct hierarchies of structured regular grids to form adaptive mesh refinement (AMR) hierarchies. Let $L \geq 1$ and $\Omega_1 \equiv \Omega \supset \Omega_2 \supset
\cdots \Omega_L$ be a nested set of subdomains of the computational
domain $\Omega$.  For simplicity, assume that each $\Omega_\ell, \; 2
\leq \ell \leq L$ is a union of non-overlapping logically rectangular regions;
these are the subregions of $\Omega$ where additional resolution is
desired.  A composite structured AMR (SAMR) grid $\CG{h}$ on $\Omega$ is a
nested hierarchy of grids $\Omega^{h_1}_1 \supset \Omega^{h_2}_2 \supset \cdots \supset \Omega^{h_{L}}_L$ consisting of $L$ levels,
with mesh spacing $h_1 > h_2 > \cdots > h_{L}$, with the coarsest grid
$\Omega^{h_1}_1$ covering $\Omega$.  Each level $\Omega^{h_\ell}_\ell$
consists of a union of non-overlapping regions, or
\textsl{patches}, at the same resolution $h_\ell$.  When there is no
risk of confusion we will drop the $\ell$ subscript and simply refer
to $\Omega^{h}_{\ell}$.  This hierarchical representation allows
operations on $\CG{h}$ to be implemented as operations on individual
levels $\Omega^{h}_{\ell}$, which in turn are decomposed into operations
on individual patches.
This property facilitates reuse of software written for regular grids.
Figure~\ref{fig:hierarchy} shows a SAMR grid with $L=3$ and one
patch on each of the refinement levels.  Note that while each
level is nested in the next coarser level, there is no requirement
that a patch at one refinement level is nested fully in a patch at
another refinement level, i.e., a fine patch at refinement level $l$
may lie over one or more coarser patches at refinement level $(l-1)$.

Having described the decomposition of a structured AMR hierarchy we now describe the discretization of (\ref{eqn:E_td})-(\ref{eqn:T_td}) on the SAMR grid hierarchy. We begin by describing the spatial discretization on a single regular grid followed by the necessary modifications for SAMR. A method of lines (MOL) approach is used where we first discretize in space to obtain a set of coupled ODEs for the variables at each spatial location, followed by discretization in time.  Only the spatial discretization for (\ref{eqn:E_td}) will be described in detail noting that a similar procedure is used to discretize (\ref{eqn:T_td}).
Integrating (\ref{eqn:E_td}) over a cell volume $\Omega_{i,j,k}$ and using the Gauss theorem for the diffusion terms we obtain:
\begin{equation}
\int\limits_{\Omega_{i,j,k}} \! \left[ \frac{\partial E}{\partial t}-\sigma_a(T^4-E)\right]\mathrm{d}V -\int\limits_{\partial \Omega_{i,j,k}} \!  (D_E \nabla E)\cdot\mathbf{n}\mathrm{d}A = 0 \label{eqn:E_int}
\end{equation}
Let $E_{i,j,k}$ and $T_{i,j,k}$ denote discrete variables collocated at cell centers indexed by $(i,j,k)$ approximating $E$ and $T$. The volumetric integral terms in (\ref{eqn:E_int}) are approximated by 
\begin{equation*}
\int\limits_{\Omega_{i,j,k}} \! \left[  \frac{\partial E}{\partial t}-\sigma_a(T^4-E)\right]\mathrm{d}V \approx \left [ \frac{\partial E_{i,j,k}}{\partial t}-\sigma_a(T^4_{i,j,k}-E_{i,j,k})\right ] h_xh_yh_z.
\end{equation*}
The surface integral diffusive flux term in  (\ref{eqn:E_int}) is evaluated approximately as
\begin{equation*}
\int\limits_{\partial \Omega_{i,j,k}} \! (D_E \nabla E)\cdot\mathbf{n}\mathrm{d}A = F_{i+\frac{1}{2},j,k}-F_{i-\frac{1}{2},j,k} + F_{i,j+\frac{1}{2},k}-F_{i,j-\frac{1}{2},k}+F_{i,j,k+\frac{1}{2}}-F_{i,j,k-\frac{1}{2}}
\end{equation*}
with
\begin{eqnarray}
F_{i-\frac{1}{2},j,k} = \int\limits_{\partial \Omega_{i,j,k}} \!  (D_E E_x)_{i-\frac{1}{2},j,k}\mathrm{d}A &\approx& (D_E)_{i-\frac{1}{2},j,k}\frac{(E_{i,j,k}-E_{i-1,j,k})}{\dx}\dy\dz \label{eqn:diffusion_coeffs1}\\
F_{i+\frac{1}{2},j,k} = \int\limits_{\partial \Omega_{i,j,k}} \!  (D_E E_x)_{i+\frac{1}{2},j,k}\mathrm{d}A &\approx&  (D_E)_{i+\frac{1}{2},j,k}\frac{(E_{i+1,j,k}-E_{i,j,k})}{\dx}\dy\dz \\
F_{i, j-\frac{1}{2},k} = \int\limits_{\partial \Omega_{i,j,k}} \!  (D_E E_y)_{i,j-\frac{1}{2},k}\mathrm{d}A &\approx& (D_E)_{i,j-\frac{1}{2},k}\frac{(E_{i,j,k}-E_{i,j-1,k})}{\dy}\dx\dz \\
F_{i,j+\frac{1}{2},k} = \int\limits_{\partial \Omega_{i,j,k}} \!  (D_E E_y)_{i,j+\frac{1}{2},k}\mathrm{d}A &\approx& (D_E)_{i,j+\frac{1}{2},k}\frac{(E_{i,j+1,k}-E_{i,j,k})}{\dy}\dx\dz \\
F_{i, j,k-\frac{1}{2}} = \int\limits_{\partial \Omega_{i,j,k}} \!  (D_E E_z)_{i,j,k-\frac{1}{2}}\mathrm{d}A &\approx& (D_E)_{i,j,k-\frac{1}{2}}\frac{(E_{i,j,k}-E_{i,j,k-1})}{\dz}\dx\dy \\
 F_{i,j,k+\frac{1}{2}} = \int\limits_{\partial \Omega_{i,j,k}} \!  (D_E E_z)_{i,j,k+\frac{1}{2}}\mathrm{d}A &\approx&(D_E)_{i,j,k+\frac{1}{2}}\frac{(E_{i,j,k+1}-E_{i,j,k})}{\dz}\dx\dy.  \label{eqn:diffusion_coeffs6}
\end{eqnarray}

The face centered diffusion coefficients, $D_E$ in
(\ref{eqn:diffusion_coeffs1})-(\ref{eqn:diffusion_coeffs6}) are computed following the description in \cite{OlsonMorel99}.
First, a face centered $T$ value based on arithmetic averaging of adjacent cell values is computed 
\begin{equation*}
T_{i-\frac{1}{2},j,k}  = \frac{1}{2}\left( T_{i,j,k}+T_{i-1,j,k} \right).
\end{equation*}
Alternatives for computing the face centered temperature described in \cite{OlsonMorel99} include geometric or parametrized arithmetic-geometric averages.
Flux matching for energy conservation with photon absorption cross sections evaluated at the face centered temperature leads to the following expression for a
harmonically averaged left face centered diffusion coefficient which is at first not flux limited:
\begin{equation*}
(D_r)_{i-\frac{1}{2},j,k} = \frac{T_{i-\frac{1}{2},j,k}^3}{3\left(z_{i,j,k}^3+z_{i-1,j,k}^3\right)}
\end{equation*}
A flux limited $D_E$ is then obtained from $D_r$ as:
\begin{equation*}
(D_E)_{i-\frac{1}{2},j,k} = \frac{2(D_r)_{i-\frac{1}{2},j,k}}{1+(D_r)_{i-\frac{1}{2},j,k}\left( \frac{|E_{i,j,k}-E_{i-1,j,k}|}{0.5h_x(E_{i,j,k}+E_{i-1,j,k})}\right )}
\end{equation*}
Similar expressions apply for  $D_E$ at other faces. We note that several alternate problem dependent choices for discretizing the diffusion coefficients are detailed in \cite{OlsonMorel99}.

At physical boundaries the values of
$E$ and $T$ are extrapolated to ghost cells using a first
order scheme and the ghost values are used in evaluating the flux
terms required at the physical boundary faces. 

% describe discretization modifications for SAMR
\subsubsection{Modifications for Structured Adaptive Mesh Refinement}
\label{sec::amr}

{\it Coarse-fine stencils:} The standard finite volume discretization of (\ref{eqn:E_td}) and (\ref{eqn:T_td}) leads to regular stencil patterns for cells that lie in the interior of a patch with each cell being connected to its six immediate neighbors normal to the faces of the cell. This promotes regular array access, minimizes cache misses and allows for the reuse of software written for regular uniform grids. However, 3D structured AMR with variable size patches on all levels can lead to complex geometric interactions between patches. Cells on the surface of a patch can lie on either the faces, edges, or corners of the patch and can be adjacent to surface cells from patches on the same and/or the adjacent coarser refinement level leading to irregular stencil patterns. One commonly used approach to restore regular array access and implicitly account for the irregular stencils is to interpolate coarse cell data into ghost cells aligned with the fine surface cells. This does require accounting for all the possible coarse-fine and fine-fine patch interactions and is programming intensive.

In Figure \ref{CoarseFineBoundaryFragments}
of the appendix, we classify the different
types of ghost cells for a representative example patch,
referring to them collectively as coarse fine boundary fragments (CFBFs).
%If there are two patches, then the CFBF belongs to the patch on the right
%(see the concave and sibling fragments cases).
As mentioned earlier, in 3D, very complex  configurations of patches can show up
and a patch can possibly have all the types of CFBFs listed.
A representative example of the types of configurations encountered is shown in Figure \ref{CoarseFineBoundaryExample}.
%As we can see from the figures,
%the CFBFs are already getting complicated,
%and the situation could get even more complex
%if there are hundreds or thousands of patches in the same level.

\begin{figure}[!ht]
\centering
\includegraphics{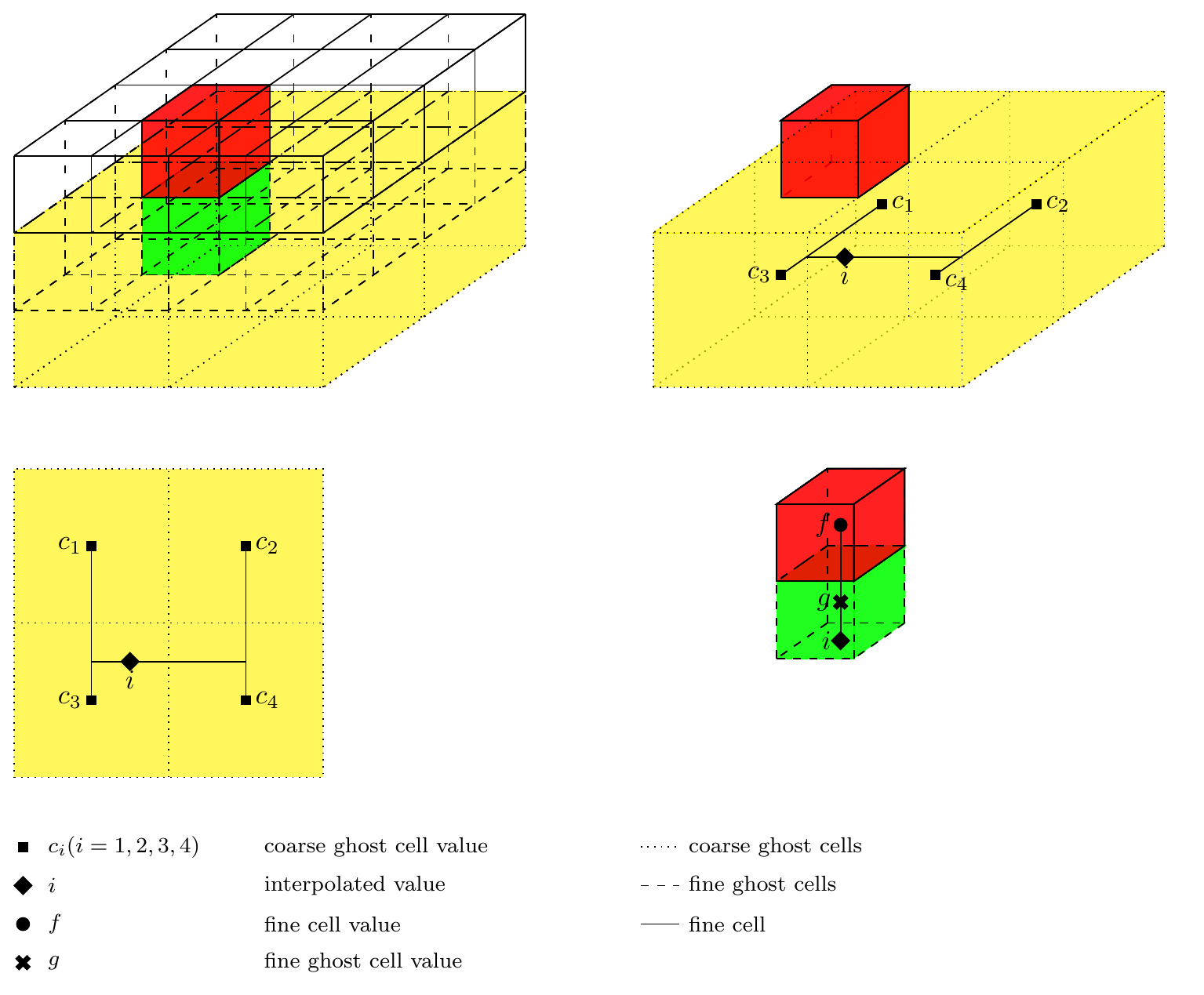}
\caption{Coarse fine boundary interpolation.}
\label{fig:interp}
\end{figure}

For the purpose of this application linear interpolation at coarse-fine
boundaries was sufficient to obtain second order accuracy as will be seen from our numerical results. A further reason to use linear interpolation is that higher
order interpolation methods in this case also suffer from the potential to overshoot and produce non-physical negative
energy and temperature values. We briefly describe linear interpolation of aligned fine ghost cells on patch faces using both coarse and fine cell data, noting
only that linear interpolation for more complex configurations required modifications to this
basic procedure that needed to be handled on a case by case basis.

Figure \ref{fig:interp} illustrates linear interpolation into a fine ghost cell using both coarse and fine values. Standard bilinear interpolation
of the four coarse ghost cell values $c_i (i=1,2,3,4)$ in the upper right figure 
(whose 2D projection is in the lower left figure) is used to obtain a coarse value
$i$ aligned with the fine interior cell $f$.
This value, while aligned is not at the fine ghost cell center.
Linear interpolation normal to the face between values $i$ and
$f$ is then used to obtain the aligned and centered fine ghost cell value
$g$ (lower right figure).
These aligned and centered fine ghost cell values are now indistinguishable from interior fine cell values for the purposes of stencil operations.

\begin{figure}[!ht]
\centering
\includegraphics{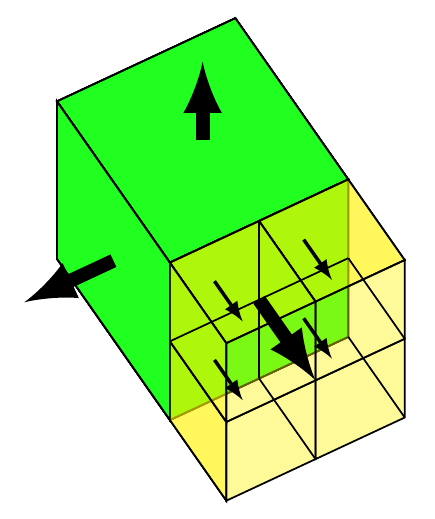}
\caption{Flux.}
\label{fig:flux}
\end{figure}

{\it Flux calculations:} Nonlinear and linear residual calculations require the calculation of the fluxes (\ref{eqn:diffusion_coeffs1})-(\ref{eqn:diffusion_coeffs6}). For flux conservation at coarse-fine boundaries, fine fluxes are averaged and are used in place of the flux calculated using coarse values only as shown in Figure (\ref{fig:flux}).

\subsection{Time discretization}
\label{sec:TimeDiscretization}
Following the spatial discretization of the previous section we obtain a system of coupled ODEs over the AMR grid that we denote by
\begin{equation}
\dot{\mathbf{u}} = \mathbf{f}(\mathbf{u}).
\label{eqn:vecform}
\end{equation}
Here $\mathbf{u}$ and $\dot{\mathbf{u}}$ respectively denote the vector of discrete unknowns and time derivatives for $E$ and $T$ across all cells. The previously cited papers, particularly \cite{mousseau-jcp-03,lowrie-jcp-04,brown-jcp-04} discuss the potential merits of various time integration schemes. In this paper we choose implicit time integration schemes due to the extreme stiffness of the systems considered as they allow us to step over the smallest time scales in the problem and instead advance at the dynamical time scale of interest. Furthermore, this enables us to leverage timestep control algorithms based on controlling accuracy from the ODE literature. We use a variable timestep backward differentiation formula of fixed order 2 (BDF2) for the numerical experiments presented. Alternate implicit time discretizations such as variable order BDF or implicit Runge-Kutta methods are possible \cite{HairerWanner96}, but this is left as a topic for future investigation. Let $\Delta t_n$ denote the $n$-th timestep and let $\alpha_n = \Delta t_n/\Delta t_{n-1}$, then a variable step BDF2 discretization of (\ref{eqn:vecform}) is given by
\begin{eqnarray}
\left ( \frac{1+2\alpha_n}{1+\alpha_n}\right ) \bunn-(1+\alpha_n)\bun+\left( \frac{\alpha_n^2}{1+\alpha_n}\right ) \bunm = \Delta t_n  \mathbf{f}(\mathbf{u_{n+1}})
\label{eqn:bdf2}
\end{eqnarray}
where $\bun$ denotes the computed solution at the $n$-th timestep. Note that for the very first timestep the solution at two previous timesteps is not available and a BDF1 method (backward Euler)
\begin{eqnarray}
\bunn = \bun + \Delta t_n  \mathbf{f}(\mathbf{u_{n+1}})
\label{eqn:be}
\end{eqnarray}
is used instead.
\subsection{Timestep control} 
In order to balance the competing objectives of accuracy and minimizing the number of required timesteps a method of timestep control is required. Within our implementation we experimented with three different algorithmic approaches to timestep control. The first approach to timestep control is based on limiting the percentage change in energy and temperature as described in \cite{Olson-jcp-07,PernicePhilip06,Rider-JCP-1999}. While simple to implement, this can be viewed as a truncation error strategy based on first order time integrators \cite{Shampine2005} potentially limiting the timestep unnecessarily. We found this to be true in practice with difficulties in choosing the necessary heuristic parameters to provide an appropriate balance between accuracy and efficiency. Choosing too small a percentage change led to small timesteps and increased solution costs while choosing too large a change led to instabilities with negative overshoots in the energy values.

The second strategy considered was the traditional ODE timestep controller \cite{Shampine2005,Brenan1987,Hindmarsh2005} based on controlling the local error per step (EPS):
\begin{equation}
\Delta t_{n+1} = \left (\frac{\epsilon_t}{|| e_n ||}\right)^{\frac{1}{k}} \Delta t_n,
\label{eqn:odets}
\end{equation}
with $k=3$ for a second order method such as BDF2. $\epsilon_t$ is a target user selected error tolerance, $e_n$ is an estimate of the local timestep error, and $|| \cdot ||$ is a norm to be specified. We note that Brown et. al. \cite{brown-jcp-04} use this strategy in the context of non-equilibrium radiation diffusion problems on uniform grids. In a series of papers \cite{Gustafsson1988,Gustafsson1994,Gustafsson1997,Soderlind2002,Soderlind2003,Soderlind2006}, Gustafsson, Soderlind and collaborators developed and analyzed timestep control algorithms including the standard ODE timestep controller described above using systematic control theory approaches. They proposed new timestep controllers that address issues of efficiency, computational stability, and tolerance proportionality that can arise with the traditional ODE timestep controller as implemented in ODE time integration packages. We consider the class of proportional-integral controllers (PI controllers) \cite{Soderlind2002} given by:
\begin{equation}
\alpha_{n+1} = \left (\frac{\epsilon_t}{|| e_n ||}\right)^{k_I} \left (\frac{|| e_{n-1} ||}{|| e_{n} ||}\right)^{k_P} \alpha_n,
\label{eqn:pc47}
\end{equation}
and in particular, we choose the $PC.4.7$ controller of \cite{Soderlind2002} with $kk_I = 0.4$ and $kk_P=0.7$ with $k=3$ for the BDF2 method. 

While we are able to control the accuracy extremely well with both of the latter approaches on uniform grids we found that modifications are required in the context of dynamic AMR. Dynamic regridding requires the transfer of solution data from an existing to a new AMR hierarchy. This introduces interpolation errors that act as non-stiff transient error components affecting the computation of the local time error estimates, $e_n$ that are based on the previous time history. The timestep controllers typically respond by dramatically reducing the timestep requiring a significant number of timesteps before the timestep is back to the dynamical timescale of the problem leading overall to an inefficient simulation. Similar effects have been noted before by Petzold\cite{Petzold1988} in a report on 1D moving grid ($r-$refinement) methods for PDEs, Trompert and Verwer \cite{Trompert1991} for 2D parabolic PDEs system simulations on structured AMR grids, and Hyman et. al \cite{Hyman2003} for hybrid moving mesh- static regridding approaches. The approach advocated in \cite{Petzold1988} is to use an implicit filtered truncation error estimator while Trompert and Verwer use a first order truncation error estimator to avoid non-stiff transients affecting the error estimates. However, as pointed out previously this can unnecessarily restrict the timestep. 

\subsection{Local time error estimation}
The timestep controllers described in (\ref{eqn:odets})-(\ref{eqn:pc47}) require an estimate of the local time error, $e_n$. At the beginning of each timestep a generalized leapfrog method  \cite{GreshoSani00} given by
\begin{eqnarray}
%\bunn^p = \bun+\left(1+\frac{\dtn}{\dtnm}\right)\dtn \dot{\bu}_n -\left(\frac{\dtn}{\dtnm}\right)^2(\bun-\bunm)
\bunn^p = \bun+\left(1+\alpha_n\right)\dtn \dot{\bu}_n - \alpha_n^2(\bun-\bunm)
\label{eqn:leapfrog}
\end{eqnarray}
 is used to provide an initial guess for the solution at the current step as well as in computing an estimate of the local error for the step. Here $\dot{\mathbf{u}}_n$ is the vector of discrete time derivative unknowns at time step $n$. Following \cite{GreshoSani00}, componentwise Taylor series expansions for the exact solution $\mathbf{u}(t_{n+1})$, the predictor $\bunn^p$, and the BDF2 solution $\bunn$ can be used to derive the following expressions for the local errors for the generalized leapfrog and BDF2 methods given by 

\noindent {\bf BDF2:}
\begin{equation}
\bunn-\bu (t_{n+1}) \approx \frac{(\dtn+\dtnm)^2}{\dtn\left (2\dtn+\dtnm \right )}\frac{\dtn^3}{6}\bun^{(3)};
\end{equation}
\noindent {\bf Generalized leapfrog:}
\begin{equation}
\bunn^p-\bu (t_{n+1}) \approx - \left(1+\frac{\dtnm}{\dtn} \right)\frac{\dtn^3}{6}\bun^{(3)}.
\end{equation}
From these expressions an approximation for the local error can be derived given by
\begin{equation}
e_{n+1} \equiv \bu_{n+1}-\bu(t_{n+1})  \approx \left( \frac{\alpha_{n}+1}{3\alpha_{n}+2}\right ) \left( \bu_{n+1}-\bu_{n+1}^p\right).\\
\label{eqn:bdf2le}
\end{equation}
In (\ref{eqn:odets}) and (\ref{eqn:pc47}) the norm of the error is a scaled max norm
\begin{equation}
||e_n||_{\infty} = \max(||e_{n}^E||, ||e_{n}^T||)
\end{equation}
where
\begin{equation}
||e_{n}^E|| = \max \left | \frac{e_{n,i,j,k}^{E}}{(E_{n,i,j,k}+\eta_E)} \right |
\end{equation}
with $e_{n,i,j,k}^E$ and $E_{n,i,j,k}$ the  the local time error estimate (\ref{eqn:bdf2le}) and the solution value of the energy density respectively for grid cell $(i,j,k)$, and $\eta_E$ a constant scaling. A similar expression is used to compute $||e_{n}^T||$. The choice of the scaled max norm instead of an $L2$ norm is important and dictated by the fact that the errors in our AMR simulations are typically local in nature when care is taken with spatial discretization. \\

\noindent{\bf Remark 1:} There appears to be an error in (3.16-250) of   \cite{GreshoSani00} leading to a different expression there for the local error estimate.\\
{\bf Remark 2:} (\ref{eqn:leapfrog}) requires an estimate for the time derivative which is computed using  the left hand side of (\ref{eqn:bdf2}). As noted in \cite{GreshoSani00} using the function evaluation directly as the expression for the time derivative term in the predictor does introduce round off errors that can lead to predictions for the timestep that are an order of magnitude lower. Furthermore, evaluation of the nonlinear function can be potentially very costly. A further reason not to use the function evaluation directly is that after regridding solution data interpolated from the old grid hierarchy is not a solution on the new grid hierarchy and any function evaluation (particularly nonlinear) will suffer significant errors. 

\section{Nonlinear solution method}
\label{sec::solver}

Efficient simulation of time dependent nonlinear multi-physics systems using implicit time integration requires attempting to minimize the total number of required timesteps over the simulation time interval through timestep control, as well as efficient solution of the nonlinear systems that arise at each timestep. Following \cite{PernicePhilip06} an inexact Newton approach, in particular, a preconditioned Jacobian Free Newton-Krylov (JFNK) method is used to solve the nonlinear systems at each timestep. JFNK methods have demonstrated their effectiveness in many applications \cite{KnollKeyes04} and in what follows we briefly summarize the general approach and then detail how our approach exploits the multilevel nature of the AMR grid hierarchies for efficiency.

\subsection{Jacobian Free Newton-Krylov Methods}

\label{sec::jfnk}

Let $\bbf:\mathbb{R}^{m}\rightarrow\mathbb{R}^{m}$ denote the nonlinear
system at each timestep and consider calculating the solution
$\bu^{\star}\in\mathbb{R}^{m}$ of the system of nonlinear equations
\begin{equation} 
\bbf(\bu^{\star})\equiv \left ( \frac{1+2\alpha_n}{1+\alpha_n}\right ) \bu^{\star}-(1+\alpha_n)\bun+\left( \frac{\alpha_n^2}{1+\alpha_n}\right ) \bunm - \Delta t_n  \mathbf{f}(\mathbf{u^{\star}})
=0\label{eqn::nonlinearsystem}
\end{equation}
that arises from (\ref{eqn:bdf2}). Classical Newton's method for solving \refeq{eqn::nonlinearsystem}
generates a sequence of approximations ${\bu^{k}}$ to $\bu^{\star}$, where
$\bu^{k+1}=\bu^{k}+\bv^{k}$ and the Newton step $\bv^{k}$ is the solution to
the system of linear equations 
\begin{equation}
\bbf'(\bu^{k})\bv^{k}=-\bbf(\bu^{k}),\label{eqn::Newtoneqns}
\end{equation} 
where $\bbf'$ is the Jacobian of $\bbf$ evaluated at $\bu^{k}$. Newton's method is
attractive because of its fast local convergence properties, but for
large-scale problems, it is impractical to solve
\refeq{eqn::Newtoneqns} with a direct method. Furthermore, it is often
unnecessary to solve (\ref{eqn::Newtoneqns}) using a tight convergence
tolerance when $\bu^{k}$ is far from $\bu^{\star}$, since the
linearization that leads to \refeq{eqn::Newtoneqns} may be a poor
approximation to $\bbf(\bu)$. Generally, it is much more efficient to
employ so-called inexact Newton methods
\cite{DemboEisentatSteihaug82}, in which the linear tolerance for
(\ref{eqn::Newtoneqns}) is selected adaptively by requiring that
$\bv^{k}$ only satisfy: 
\begin{equation}
\norm{\bbf(\bu^{k})+\bbf'(\bu^{k})\bv^{k}}\leq\eta_{k}\norm{\bbf(\bu^{k})}
\label{eqn::inexactNewtoneqns}
\end{equation}
for some $\eta_{k}\in(0,1)$ \cite{DemboEisentatSteihaug82}. Appropriate choice of the
forcing term $\eta_{k}$ can lead to superlinear and even
quadratic convergence of the iteration  \cite{EisenstatWalker94}. The algorithm we use can be found in \cite{KelleyIterativeBook}. 

\noindent{\bf Remark:} A modification introduced into the general inexact Newton algorithm above due to the application was to constrain the search direction vectors $\bv^k$ so that the update vectors $\bu^k+\bv^k$ maintain the positivity of $E$ and $T$.

While any iterative method can be used to find a $\bv^{k}$ that satisfies
\refeq{eqn::inexactNewtoneqns}, Krylov subspace methods are distinguished
by the fact that they require only matrix-vector products to proceed.
These matrix-vector products can be approximated by a
finite-difference version of the directional (G\^ateaux) derivative as:
\begin{equation} 
\bbf'(\bu^{k})\bv\approx\frac{\bbf(\bu^{k}+\varepsilon \bv)-\bbf(\bu^{k})}{\varepsilon},\label{eqn:Jac-free}
\end{equation} 
which is especially advantageous when $\bbf'$ is difficult to compute or expensive
to store (both being the case here due to the presence of
multiple grid patches across the AMR grid hierarchy). 

Several approaches exist to compute the differencing parameter $\varepsilon$ \cite{KelleyIterativeBook}. In general, the main consideration in the choice of $\varepsilon$ is the accuracy of the approximation (\ref{eqn:Jac-free}). A further consideration that cannot be relaxed in this application is maintaining the positivity of the component $E$ and $T$ values for the vector $\bu^k+\varepsilon \bv$ in (\ref{eqn:Jac-free}). The conditions when the positivity constraint could be violated seem to occur most frequently immediately after regridding when interpolation and coarsening of solution data between old and new AMR grid hierarchies introduces transient numerical errors. The most robust choice (while more expensive than some other alternatives) in our numerical experiments was to compute 
\[
\varepsilon = \begin{cases}\frac{\sqrt{\epsilon_{\mathrm{mach}}}\langle \bu,\bv \rangle}{\|\bv\|^2}  & \mbox{if} \langle \bu,\bv \rangle > b\bu_{min}\|\bv\|_1 \\
                                                 \frac{\sqrt{\epsilon_{\mathrm{mach}}}\bu_{min} sign( \langle \bu,\bv \rangle ) \|\bv\|_1}{\|\bv\|^2} & \mbox{otherwise},
                         \end{cases}
\] 
where $\epsilon_{\mathrm{mach}}$ is machine precision, $\bu_{min}$ is set to $10^{-6}$, $\|\cdot\|$ and $\|\cdot\|_1$ 
are the discrete $L2$- and $1$-norm respectively. In our numerical experiments performed with
double precision arithmetic, $\varepsilon$ is typically on the order of $10^{-8}$.

Among the Krylov methods appropriate for non-symmetric definite systems we choose the GMRES  method for its robustness \cite{saad-gmres}. A further advantage of GMRES when
employed as part of a JFNK method is that the Krylov vectors  are normalized, $\|\bv\|=1$, bounding the error introduced in the difference approximation of (\ref{eqn:Jac-free})  whose leading error term is
proportional to $\varepsilon\| \bv\|^{2}$ \cite{knoll-aiaa-ns}. The main drawbacks commonly cited with GMRES are that it can potentially require the storage of a large number of Krylov vectors making it memory intensive and that it can be compute intensive as its computational complexity is proportional to the square of the number of GMRES iterations at each Newton step. We limit the potential impact on efficiency of these characteristics of GMRES by minimizing the number of GMRES iterations required at each Newton step through a combination of inexact Newton methods (to prevent oversolving) and strong physics based preconditioning to improve the conditioning of the systems. The next section focuses on the details of the preconditioner and its application. 

\subsection{Preconditioners}
JFNK allows us to focus on developing effective preconditioners. Right-preconditioning of the Newton equations is used, i.e., we solve
\begin{displaymath}
(\bbf'(\bu^k)P^{-1})P \bv^k = -\bbf(\bu^k).
\end{displaymath}
where $P$ is the preconditioner. From (\ref{eqn:Jac-free}), this requires the Jacobian-vector products
\begin{displaymath}
\bbf'(\bu^k)(P^{-1}\mathbf{w}) \approx \frac {\bbf(\bu^k+\varepsilon P^{-1}\mathbf{w})-\bbf(\bu^k)}{\varepsilon}
\end{displaymath}
within GMRES.
These are computed in two steps. An application of the preconditioner, $\mathbf{y}=P^{-1}\mathbf{w}$, yields the vector $\mathbf{y}$ that is then used to compute
$\frac {\bbf(\bu^k+\varepsilon \mathbf{y})-\bbf(\bu^k)}{\varepsilon}$. Having described the general approach we now describe how the preconditioner $P$ is constructed. For preconditioning  the Jacobian systems at each Newton step are approximated by
\renewcommand{\arraystretch}{0.5}
\begin{displaymath}
\mathcal{L}^k=
\left(\begin{array}{cc}
       (1+\sigma_a \beta_n)I - \beta_n \div D_E^{k}\grad  & -\sigma_a \beta_n(T^k)^3 I\\
       -\sigma_a\beta_n I & (1+\sigma_a\beta_n(T^k)^3)I -\beta_n\div D_T^k\grad 
       \end{array}\right)\label{eqn:operator}
\end{displaymath}
where $D_E^{k}$ and $D_T^{k}$ are the diffusion coefficients frozen at the previous Newton iterate values and $\beta_n= \left(\frac{\alpha_n+1}{2\alpha_n+1}\right) \Delta t_n$.
Several options for preconditioning are described in the literature. We choose to extend the following multiplicative splitting described in \cite{Mousseau00} to AMR grids. 
\begin{displaymath}
\mathcal{L}^k\approx \mathcal{P}_1^k\mathcal{P}_2^k
\end{displaymath}
where
\begin{displaymath}
\mathcal{P}_1^k=
\left(\begin{array}{cc}
       I - \beta_n \div D_E^{k}\grad & 0\\
       0 & I -\beta_n\div D_T^k\grad 
       \end{array}\right)
\end{displaymath}
and
\begin{displaymath}
\mathcal{P}_2^k=
\left(\begin{array}{cc}
       (1 +\sigma _a\beta_n)I  & -\sigma_a\beta_n(T^k)^3 I\\
       -\sigma_a\beta_n I & (1 +\sigma_a\beta_n(T^k)^3 )I
       \end{array}\right)
\end{displaymath}

Systems involving $\mathcal{P}_2^k$ only require local cell by cell inversion of block $2\times 2$ systems which is an embarrassingly parallel operation over the whole AMR grid. Systems involving $ \mathcal{P}_1^k$ on the other hand contain two diagonal variable coefficient elliptic operators. Such systems could be approximately solved using a variety of methods such as block Jacobi, ILU variants, or multigrid on a per patch or per refinement level basis. However, methods that fail to account for inter-level couplings cannot eliminate any global error modes that span the refinement levels resulting in preconditioner performance that progressively degrades as the number of refinement levels is increased. In principle an algebraic multigrid or semistructured multigrid method could be used over the whole AMR grid hierarchy to address this. Such methods are attractive particularly on unstructured grids where defining geometric grid coarsenings is hard. However, in the structured AMR context  such methods would require the formation of irregular stencil operators across refinement levels, a task that in practice is extremely programming intensive and error prone for finite volume discretizations with coarse-fine boundary interpolation. Furthermore, such methods in general ignore the multilevel structure of the AMR grid that already exists. The multilevel method we choose to invert the components of $\mathcal{P}_1^k$ with is the Fast Adaptive Composite Grid (FAC) method \cite{McCormickThomas86,McCormick89} of McCormick et. al. that extends techniques from multigrid on uniform grids to exploit the natural multilevel structure of SAMR grid hierarchies. 

\subsection{Fast Adaptive Composite Grid (FAC) method}
FAC solves problems on AMR grids by
combining smoothing on refinement levels with a coarse grid solve
using an approximate solver, such as a V-cycle of multigrid. In order to describe the algorithm we establish the following notation.

\begin{itemize}
  \item $R_{\ell}:\Omega_{\ell+1}^h \rightarrow \Omega_{\ell}^h$ 
        and $P_{\ell}:\Omega_{\ell-1}^h \rightarrow \Omega_{\ell}^h$ 
        respectively denote restriction and interpolation operators between
        adjacent refinement levels. Here we  use bilinear interpolation for $P_{\ell}$ and simple averaging 
        for $R_{\ell}$.

  \item $R_{\ell}^c:\Omega^c \rightarrow \Omega^h_{\ell}$ and
        $P_{\ell}^c:\Omega^h_{\ell} \rightarrow \Omega^c$ respectively 
        denote restriction and interpolation operators between the 
        composite AMR grid $\Omega^c$ and an individual refinement level $\Omega^h_{\ell}$. In practice these are expressed in terms of
        the inter-level operators $R_{\ell}$ and $P_{\ell}$.
        
  \item $\mathcal{L}^c$ represents a composite fine grid discrete operator discretized over the AMR grid hierarchy 
       $\Omega^c$ and represents for example one of the components of $\mathcal{P}_1^k$.
       $\mathcal{L}^{\ell}$ approximates $\mathcal{L}^c$ on level $\ell$.
\end{itemize}
With this notation we can specify one V($m$,$n$) sweep of the FAC Method as in
Algorithm~\ref{alg::fac}.
\begin{algorithm}[ht]
\caption{\FAC}
\label{alg::fac}
\begin{algtab*}
\algbegin
Initialize:  $r^c = f^c - \mathcal{L}^cu^c$; $f^{\ell} = R_{\ell}^c r^c $ \\
\algforeach{$\Omega^{h}_{\ell}, \; \ell = L, \ldots,  2 $}
Smooth $m$ times on: $\mathcal{L}^{\ell} e^{\ell} = f^{\ell}$                          \\
Correct\hspace{2pt}: $u^c = u^c+P_{\ell}^c e^{\ell}$                      \\
Update\hspace{2pt}: $r^c = f^c-\mathcal{L}^cu^c$                          \\
Set\hspace{20pt}: $f^{\ell-1} = R^c_{\ell-1} r^c$                         \\
\algend
Solve\hspace{10pt}: $\mathcal{L}^1 e^{1} = f^{1}$                         \\
Correct: $u^c = u^c + P_{1}^c e^{1}$                                      \\
\algforeach{$\Omega^{h}_{\ell}, \; \ell = 2, \ldots, L$}
Smooth $n$ times on: $\mathcal{L}^{\ell} e^{\ell} = f^{\ell}$                          \\
Correct\hspace{2pt}: $u^c = u^c + P_{\ell}^c e^{\ell}$
\end{algtab*}
\end{algorithm}

Algorithm~\ref{alg::fac} makes clear the \emph{multiplicative} nature
of FAC: the residual is updated with the latest correction information
before each smoothing pass can proceed.  To be fully effective, each
smoothing pass must properly account for the data dependencies among
different patches within a refinement level as well as coarse-fine data dependencies.  In our calculations, we
use red-black Gauss-Seidel smoothing on each refinement
level; we also have the capability to use damped point Jacobi or block
Gauss-Seidel smoothing.  The correction steps require synchronization
of the composite grid solution to make it consistent on all refinement
levels.  Note that the residual update can, in principle, be computed
on only the most recently corrected refinement level plus a small
border on the next coarser level, but we have found that residual
evaluation is not expensive enough to justify this optimization.  On
the coarsest level, we can use one V-cycle of a multigrid method. However, our numerical results will simply use smoothing on the coarsest level.

\section{Dynamic Regridding}
Dynamic AMR requires changing the patch hierarchy as the simulation evolves in response to a changing error profile over the domain. This leads to refinement and derefinement of subdomains. Heuristic error indicators given by 
\begin{equation}
\tau^c_{i,j,k} = \frac{h_x^2|(E_{xx})_{i,j,k}|+h_y^2|(E_{yy})_{i,j,k}|+h_z^2|(E_{zz})_{i,j,k}|}{0.1 \max_{\substack{i,j,k}}|E_{i,j,k}|}
\end{equation}
and 
\begin{equation}
\tau^g_{i,j,k} = \frac{h_x|(E_{x})_{i,j,k}|+h_y|(E_{y})_{i,j,k}|+h_z|(E_{z})_{i,j,k}|}{0.1 \max_{\substack{i,j,k}}|E_{i,j,k}|}
\end{equation}
are used to identify cells with high curvature and gradient in the energy density based on a user specified threshold. Here $E_x$ and $E_{xx}$ denote the finite difference approximations to the gradient and Laplacian in the $x$ direction with a similar notation being used for the other terms. Tagged cells are grouped into patches based on the Berger-Rigoutsos algorithm \cite{BergerRigoutsos91} as implemented in the SAMRAI package \cite{samrai}. Regridding is critical in time dependent calculations to minimize the total computational cost as the simulation evolves. However, it is an expensive operation and can introduce transient error into the simulation that leads to reduced timestep sizes and an overall increase in the number of timesteps required unless handled carefully. The error indicators above are calculated at fixed intervals ( typically every 10 timesteps ) and a regrid operation is triggered only when necessary. 

Following the creation of a new grid hierarchy, interpolation of data from the old patch hierarchy to the new is required. The converged solution at a timestep on the old grid hierarchy is in general not a solution on the new grid hierarchy after interpolation as is often evidenced by the jump in the nonlinear residual. High order interpolation could in principle minimize this effect but leads to non-physical negative undershoots in the interpolated values of the energy and temperature and is non-conservative. Hence, in our calculations we use conservative linear refinement to eliminate this possibility. Following solution interpolation onto the new AMR hierarchy time integration can be restarted. A cold restart approach to regridding uses the interpolated solution as the initial condition to restarting the simulation starting with a single stage method such as backward Euler. A warm restart continues the simulation using interpolated time history vectors. We adopt a modified version of the latter approach by first performing a re-solve at the current timestep using the interpolated time history vectors as an initial guess. This however does not eliminate small discontinuities in the time derivatives and introduces non-stiff transient error components that are detected by the timestep controller leading to reduced time step sizes. Furthermore, since the timestep controllers rely on accurately estimating the local time error (\ref{eqn:bdf2le}) which depends on quantities that can only be obtained by interpolation from the old grid hierarchy the estimate  (\ref{eqn:bdf2le}) is necessarily wrong for the first timestep and is not used.

\section{Numerical results}
\label{sec:numerics}
\begin{figure}[!ht]
\centering
\includegraphics{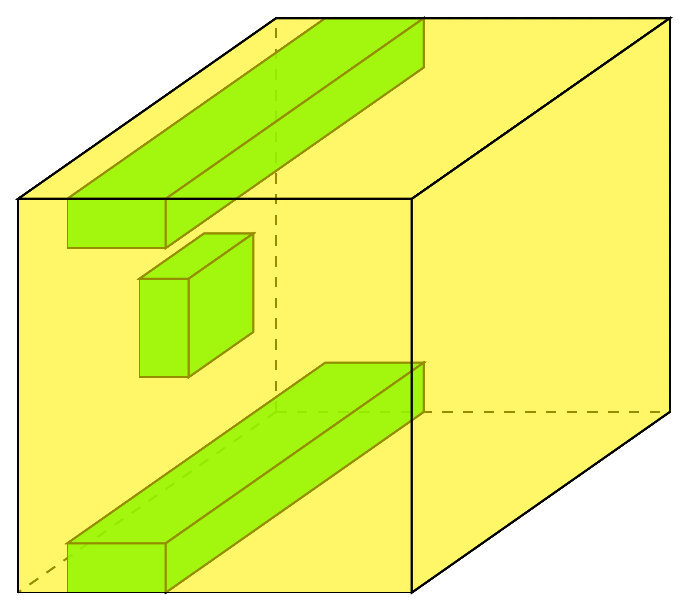}
\caption{Material configuration.}
\label{fig:material}
\end{figure}

Numerical results are presented for a carefully chosen model problem that is designed to test the performance of the implicit time integrator and its components on fairly complex 3D AMR hierarchies with challenging variations in the material properties. The domain is the unit cube with regions containing two materials with 
\begin{equation*}
z(x,y,z) = 
\begin{cases}
10, \quad (x,y,z) \in
\begin{cases}
   [0.0625, 0.2]  \times [0.375,0.625] \times [0.375,0.625], \\
   [0.125, 0.375] \times [0.0,1.0] \times [0.0,0.125], \\
   [0.125, 0.375] \times [0.0,1.0] \times [0.875,1.0],
\end{cases} \\
1, \quad \mathrm{otherwise},
\end{cases}
\end{equation*}
as illustrated in Figure \ref{fig:material}. At the initial time the energy density and material temperature are constant over the domain with $E_0 = 10^{-5}$, $T_0 = (E_0)^{\frac{1}{4}}$. On the left and right faces of the domain Robin boundary conditions (\ref{eqn:bcs}) for the energy density are imposed with $R=1$ and $R=0 $ respectively. Zero Neumann conditions are imposed for the energy density on all other faces.  For the material temperature, zero Neumann boundary conditions are imposed on all boundaries. As the simulation progresses the left boundary of the domain is heated up and a steep Marshak wave \cite{Marshak1958} front for the energy density flows from the left to the right of the domain. AMR is used to resolve this wave front accurately with the implicit time integration schemes described enabling us to timestep efficiently and accurately at the dynamical time scale of the problem. As the wave front hits the regions with high $z$ number fairly complex AMR grid configurations are generated due to the $z^3$ dependence of the diffusion coefficients. We locate these regions in our numerical experiment close to the left boundary to generate complex AMR configurations early on. 

All numerical results presented are based on simulations integrated to final time $t=1.0$ with a variable step BDF2 implicit integrator except for the $256^3$ equivalent simulations where the final time $t=0.5$ due to limits on the computational resources available to us. The user set tolerance $\epsilon_t$ in (\ref{eqn:odets}), (\ref{eqn:pc47}) is set to $5.0\times10^{-4}$ for grids with resolution equivalent to a $32^3$ uniform grid and decreased by a factor of two each time the grid resolution is doubled. All the numerical results reported use the $PC.4.7$ timestep controller described earlier. At each timestep a JFNK method is used with the forcing term, $\eta_k$ set to according to the Eisentat-Walker algorithm as described in \cite{KelleyIterativeBook}. The relative and absolute tolerances are $10^{-12}$ and $10^{-10}$ respectively in the nonlinear solver. The Krylov method used is GMRES with the maximum Krylov space dimension set to $50$ though this limit is never reached in practice. GMRES is right preconditioned with one iteration of the physics based preconditioner as described previously with single V(1,0) cycles of FAC being used to approximately invert the diffusion components for temperature and energy in the preconditioner. Red-black Gauss-Seidel is the smoother used on all refinement levels. All computations are performed in double precision on a Cray XK6 machine at ORNL. The infrastructure for SAMR was provided by the SAMRAI framework \cite{samrai}, the JFNK solver is a slight modification of the PETSc \cite{Balay_etal01,petsc,petsc-web-page} SNES implementation of the Eisenstat-Walker inexact Newton algorithm and the FAC solver was part of the SAMRSolvers package developed by the authors. 

\begin{figure}[!ht]
\begin{center}
\begin{tabular}{ccc}
\includegraphics[width=0.5\linewidth]{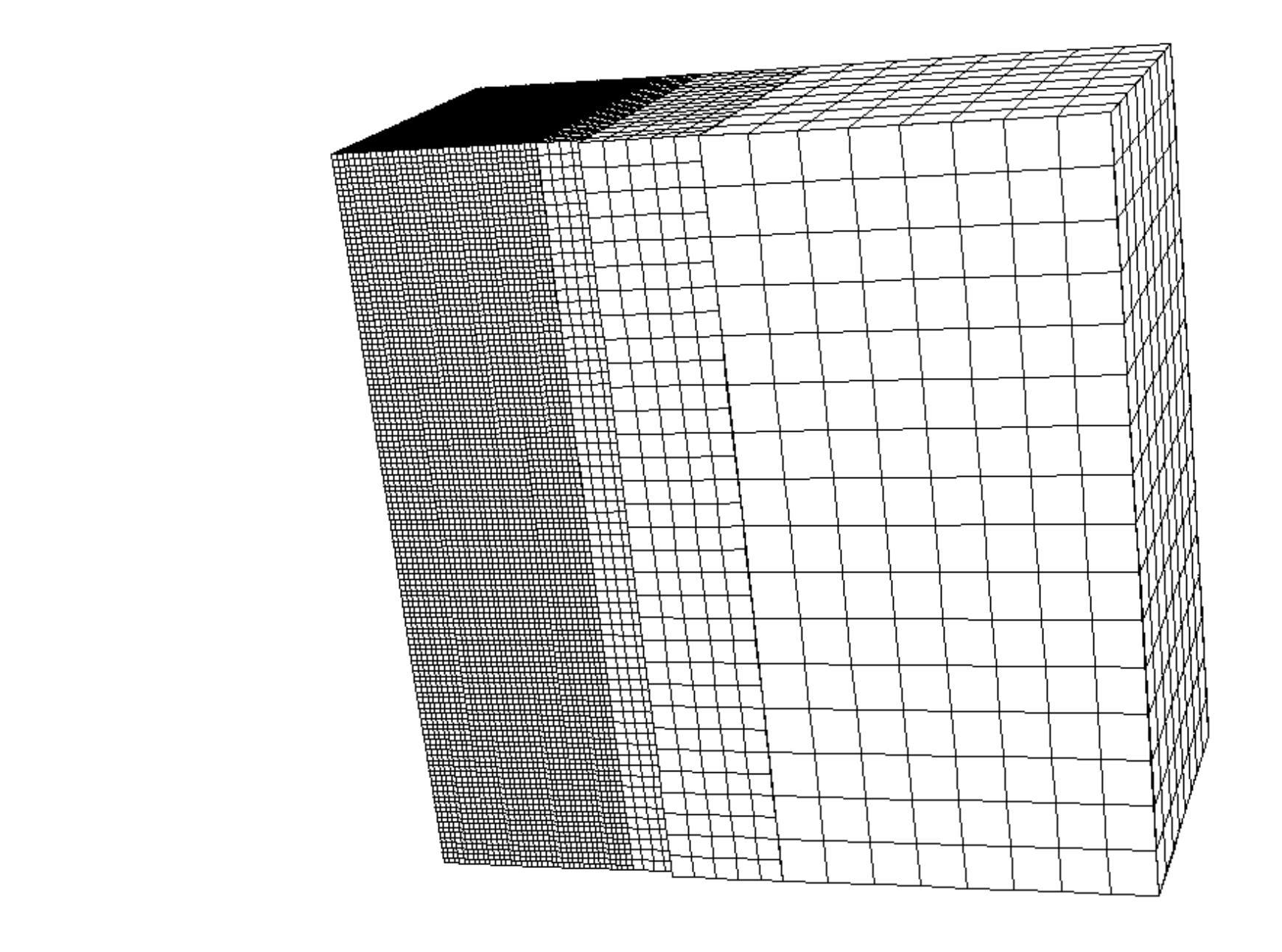} &
\includegraphics[width=0.5\linewidth]{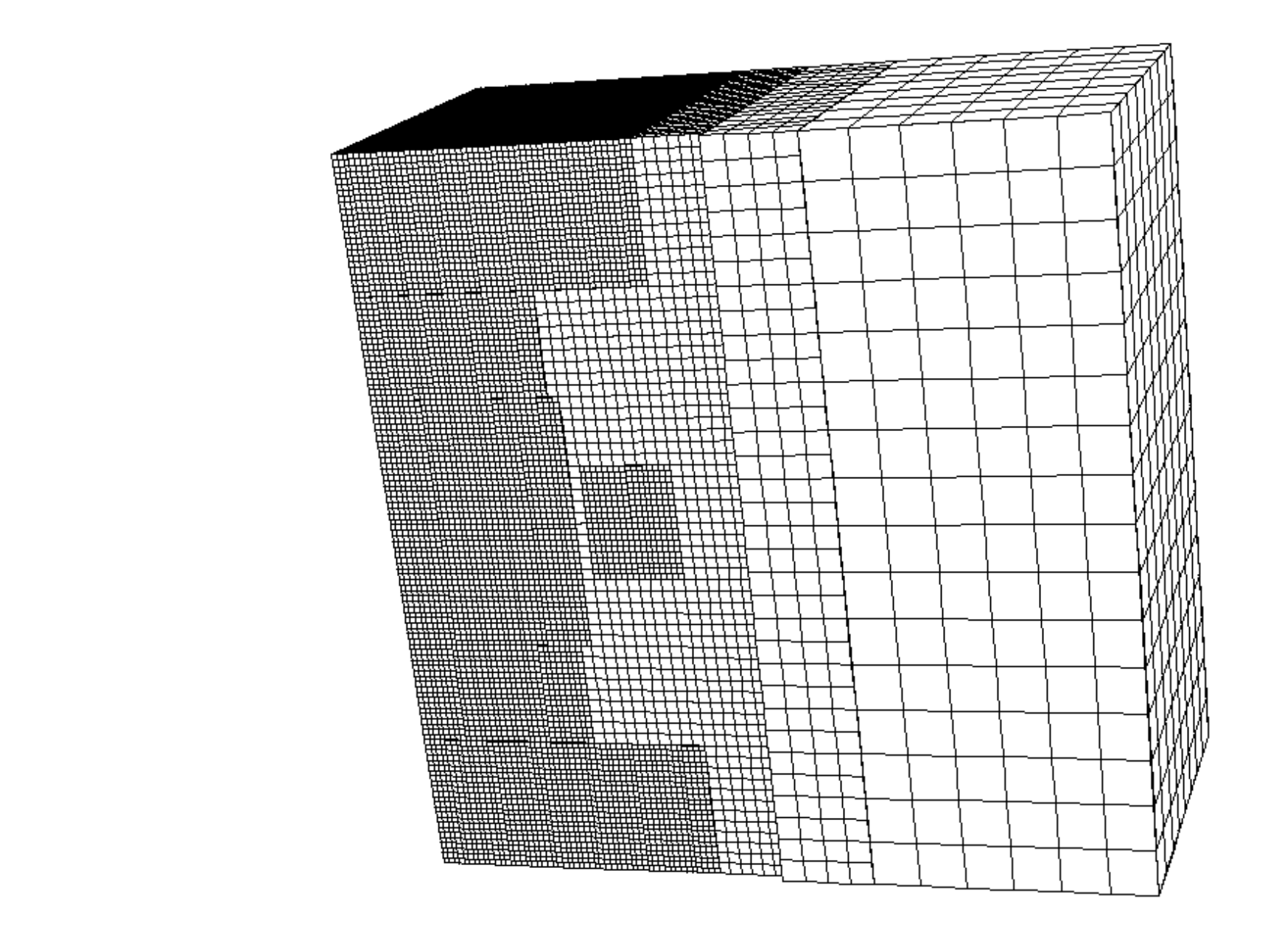} \\
\includegraphics[width=0.5\linewidth]{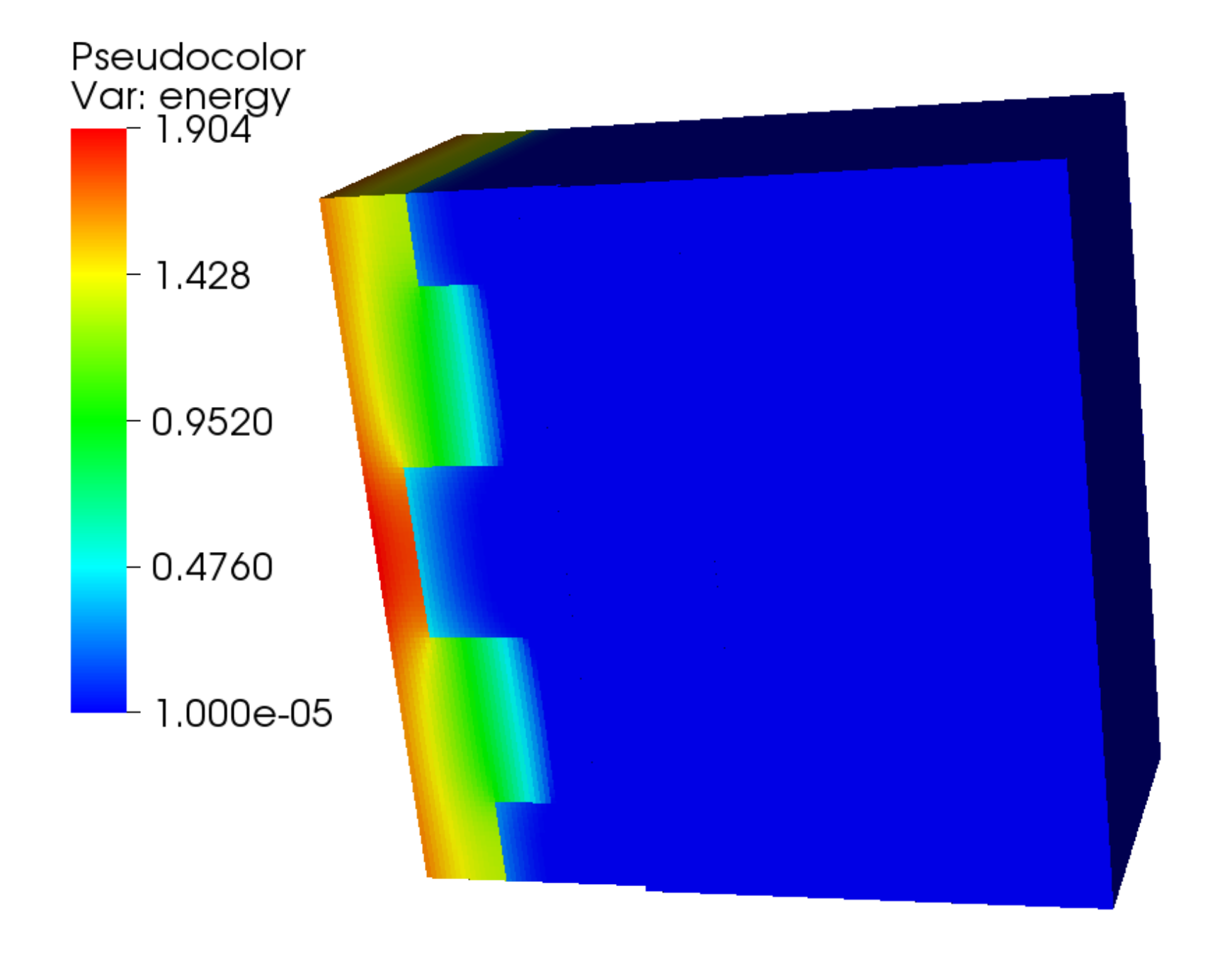} &
\includegraphics[width=0.5\linewidth]{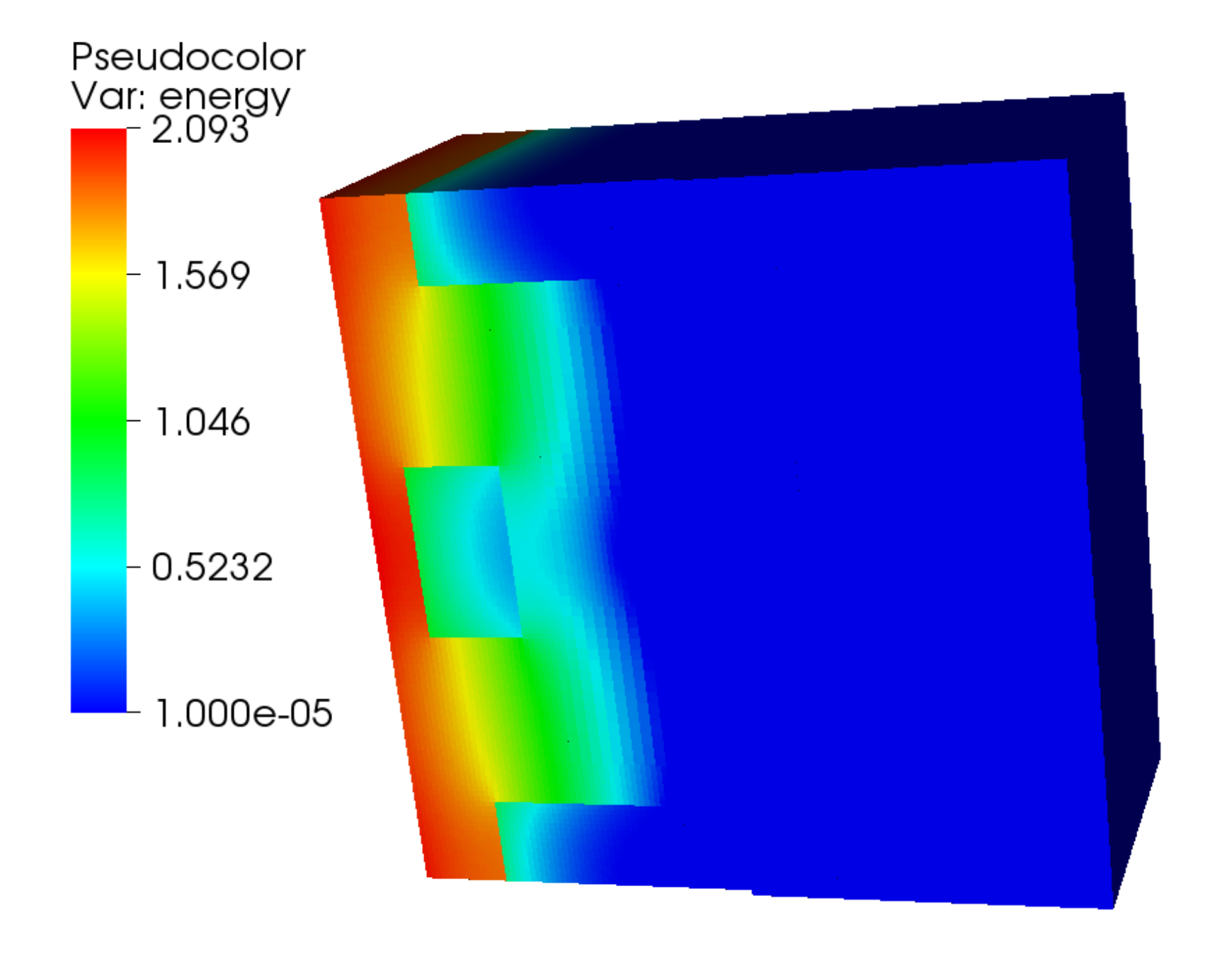} \\
\includegraphics[width=0.5\linewidth]{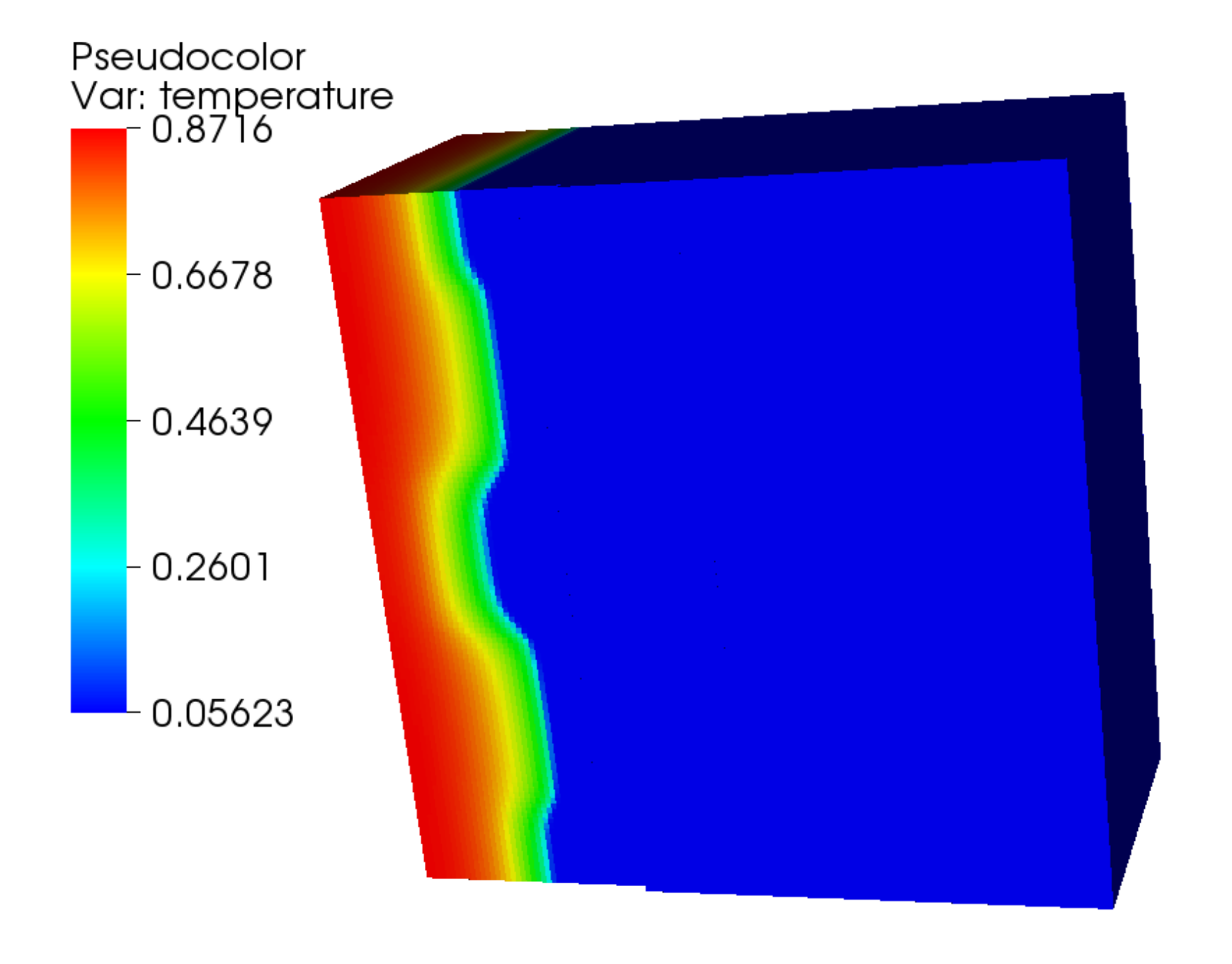} &
\includegraphics[width=0.5\linewidth]{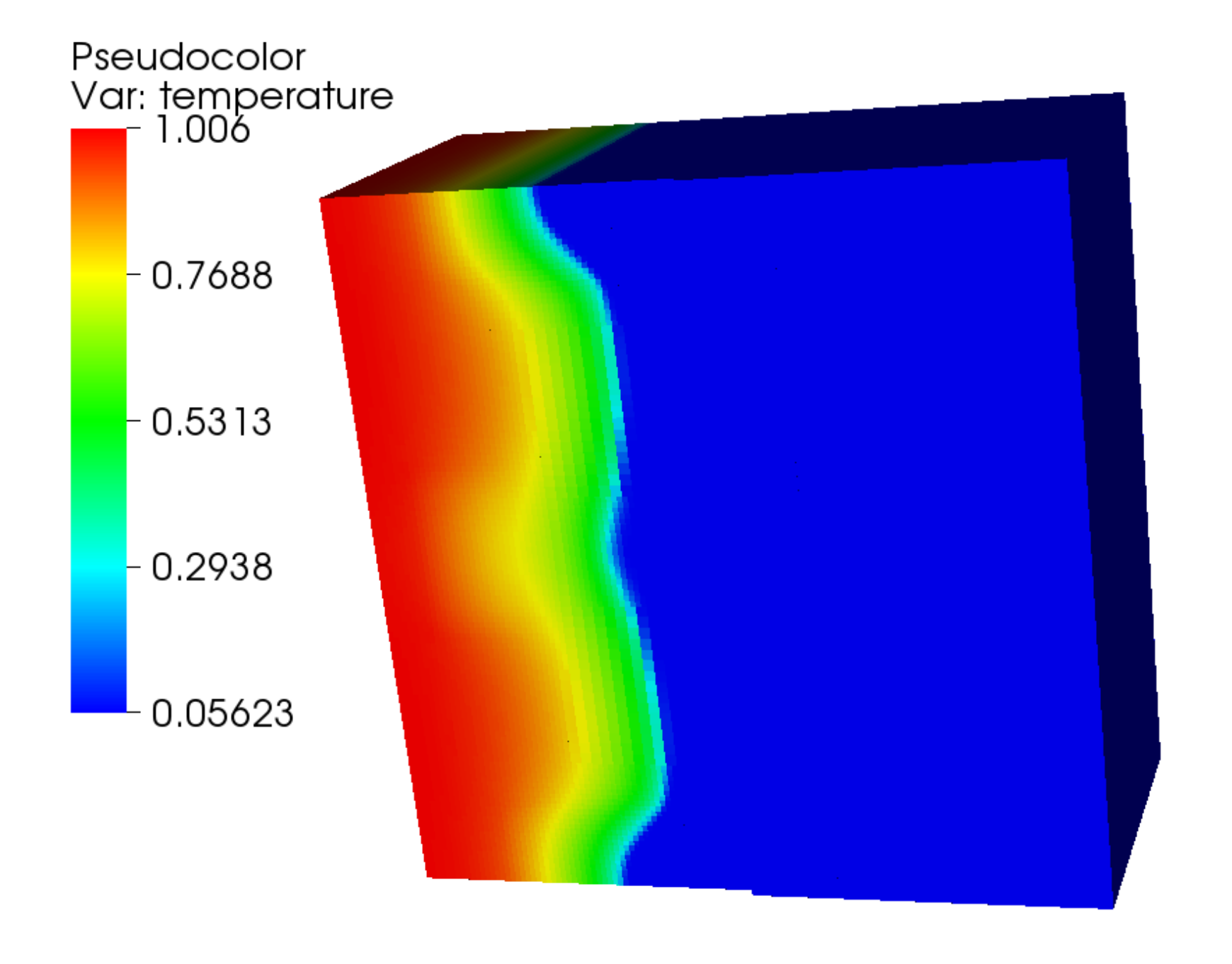} \\
t = 0.5 & t = 1.0 \\
\end{tabular}
\end{center}
\caption{Evolution of AMR mesh, energy and temperature on a 16b4l grid
on the domain $ [0,1] \times [0.52,1] \times [0,1] $.
From top to bottom, the figures represent the AMR mesh, energy density, and temperature
respectively.}
\label{fig:evolution}
\end{figure}

Snapshots of the time evolution of an AMR mesh and solutions are presented
in Figure \ref{fig:evolution}. In the figures, the coarsest level is an uniform $ 16^3 $ grid,
with the finest level resolution equivalent to that of an uniform $ 128^3 $ grid.
As the figures show, the AMR grid dynamically tracks the movement of the Marshak wave front, with
the energy plots in particular showing the relatively slow heating of the material with $ z(x,y,z) = 10 $.

\subsection{Efficiency study}

\begin{table}[h]
\parbox{0.45\linewidth}{
\centering
\renewcommand{\arraystretch}{1.3}
\begin{tabular}{|c|c|c|c|c|c|}
\hline
   Levels      & 1 & 2 & 3 & 4 & 5 \\
\hline\hline
$ 16^3     $ & - & 7.71 & 7.21 & 6.63 & 6.92 \\ \hline
$ 32^3     $ & 7.68 & 7.22 & 6.61 & 6.91 & - \\ \hline
$ 64^3     $ & 7.15 & 6.63 & 6.92 & - & -      \\ \hline
$ 128^3  $ &  6.43 & 6.86 & - &   -    &    -     \\ \hline
$ 256^3  $ &  6.72 & - & - &   -    &    -     \\ \hline
\end{tabular}
\caption{\label{case1-li} Average linear iterations}
}
\hfill
\parbox{0.45\linewidth}{
\centering
\renewcommand{\arraystretch}{1.3}
\begin{tabular}{|c|c|c|c|c|c|}
\hline
   Levels      & 1 & 2 & 3 & 4 & 5 \\
\hline\hline
$ 16^3      $ & - & 2.99 & 2.96 & 2.64 & 2.71 \\ \hline
$ 32^3      $ & 2.99 & 2.96 & 2.63 & 2.71 & -  \\ \hline
$ 64^3      $ & 2.96 & 2.64 & 2.71 & - & -      \\ \hline
$ 128^3    $ & 2.61 & 2.72 & - &   -    &    -     \\ \hline
$ 256^3  $ &  3.04 & - & - &   -    &    -     \\ \hline
\end{tabular}
\caption{\label{case1-nli}Average nonlinear iterations}
}
\end{table}%

Having described the numerical experiment setup we now study the efficiency of the implicit time integration scheme and the nonlinear and linear solvers on various AMR grids.
Table \ref{case1-li} presents the average number of linear iterations required at each timestep as we vary the number of refinement levels when using FAC as a preconditioner. Along each row we have the number of refinement levels for the AMR grid and along each column the number of grid cells on the coarse base grid. Moving diagonally in Table \ref{case1-li} from the lower left to the upper right we see that for the same effective finest grid resolution the average number of linear iterations required to solve a problem on a refined grid remains approximately constant. For example, on average 6.43 linear iterations are required per timestep to advance the solution on a uniform grid with $128^3$ grid cells and on average approximately the same number of iterations is required to advance the solution on adaptively refined grids with $16^3$,  $32^3$, $64^3$, and $128^3$  base grids with 4, 3, and 2 levels of refinement respectively. This demonstrates the level independent performance of the multilevel physics based preconditioner without which we would have seen the average number of linear iterations increase simply as a result of moving from a uniform to an AMR grid with the same resolution. Moving vertically down each column we see that increasing grid resolution marginally appears to decrease the number of linear iterations required. We infer from this that the timestep selection scheme is behaving correctly and that the FAC preconditioning is enabling the Krylov solver to deliver condition number independent convergence at each Newton step. For all our simulations a V(1,0) (slash) cycle of FAC was used as no significant improvement was  seen by increasing the number of pre- or post- smoothing steps in FAC. In Table \ref{case1-nli} we present the average number of {\it nonlinear} iterations per timestep as the base grid resolution and number of refinement levels is varied. Here again we see little variation in the required number of nonlinear iterations for the same effective fine grid resolution. 

\begin{figure}[!ht]
\begin{center}
\includegraphics[width=0.9\linewidth]{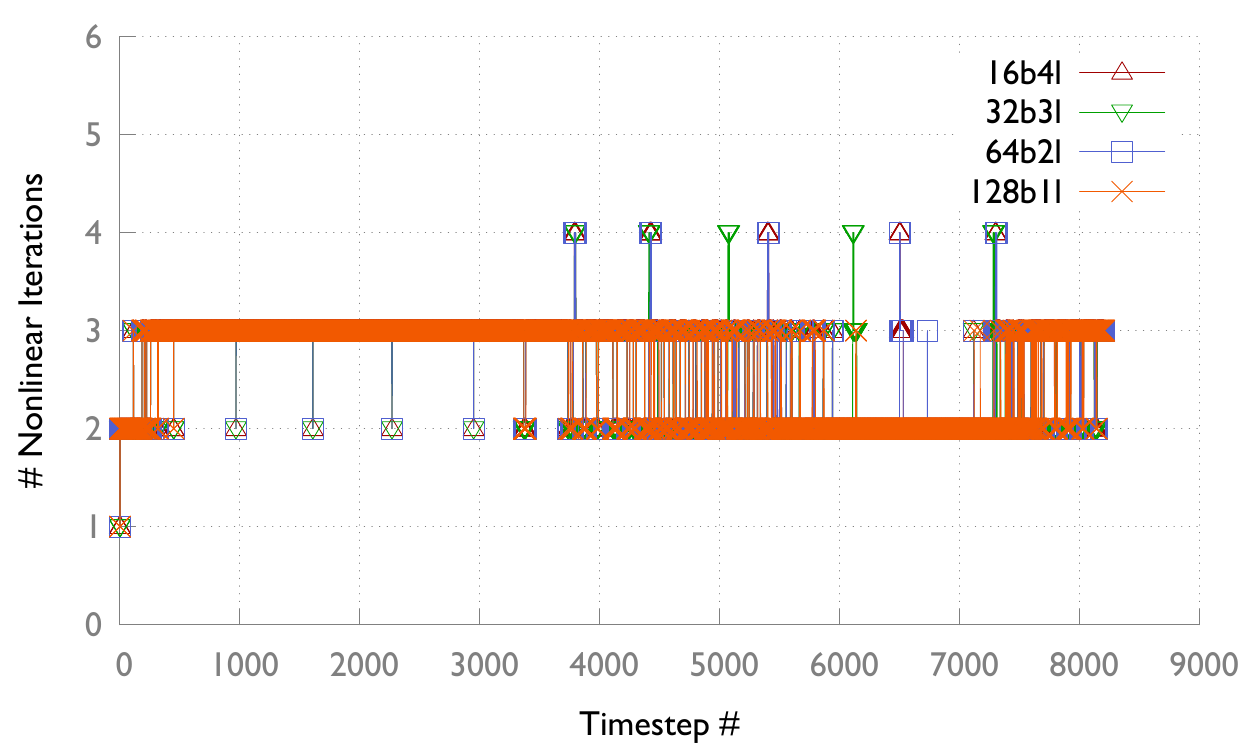}
\end{center}
\caption{Time history of nonlinear iteration counts for different AMR grids with equivalent resolution}
\label{fig:nlihistory}
\end{figure}

\begin{figure}[!h]
\begin{center}
\includegraphics[width=0.9\linewidth]{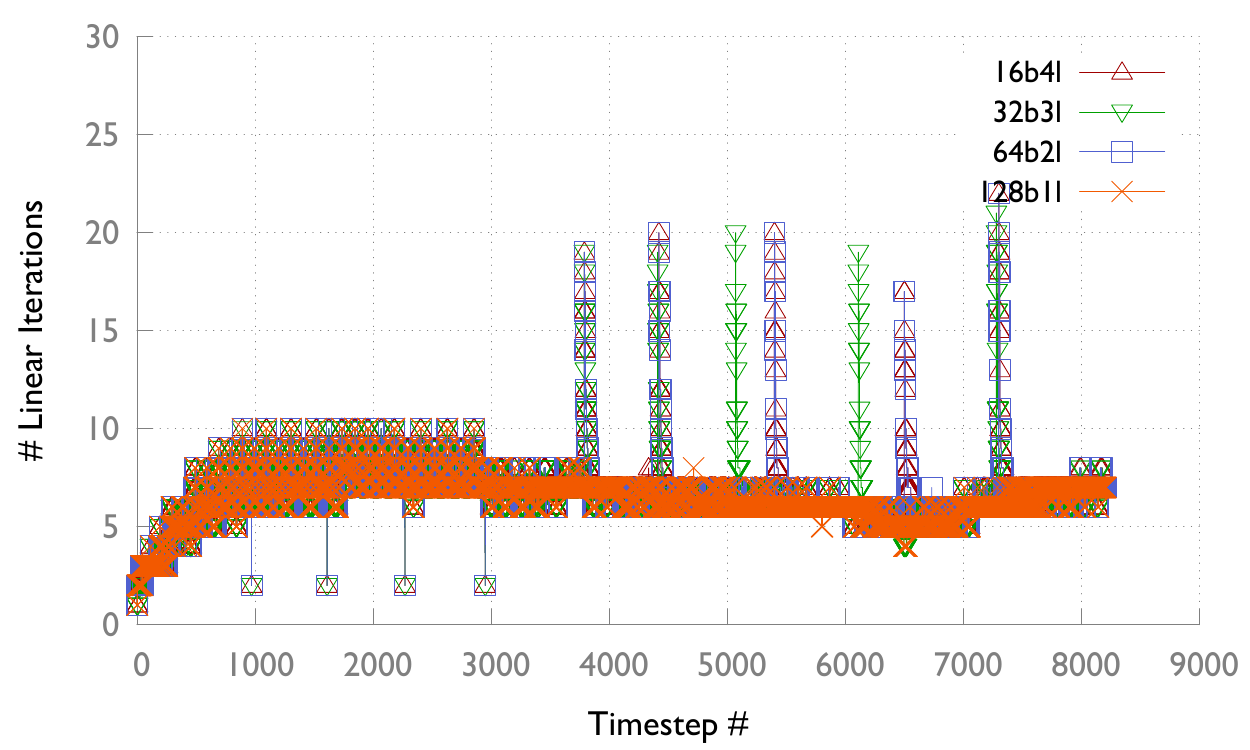} 
\end{center}
\caption{Time history of linear iteration counts for different AMR grids with equivalent resolution}
\label{fig:lihistory}
\end{figure}

In figures (\ref{fig:nlihistory}) and (\ref{fig:lihistory}) the per time step nonlinear and linear iteration counts for AMR grids with resolution equivalent to a $128^3$ uniform grid are shown. The key `16b4l' in the figures refers to an AMR grid with a coarse base grid having $16^3$ grid cells and 4 refinement levels including the base grid. The other keys are to be interpreted similarly. As can be seen the performance of the nonlinear and linear solvers is fairly comparable for the different grids. Major differences in iteration counts are only present at the regrid events where spikes are seen in the plots.

\begin{figure}[h]
\centering
\includegraphics[width=0.9\textwidth]{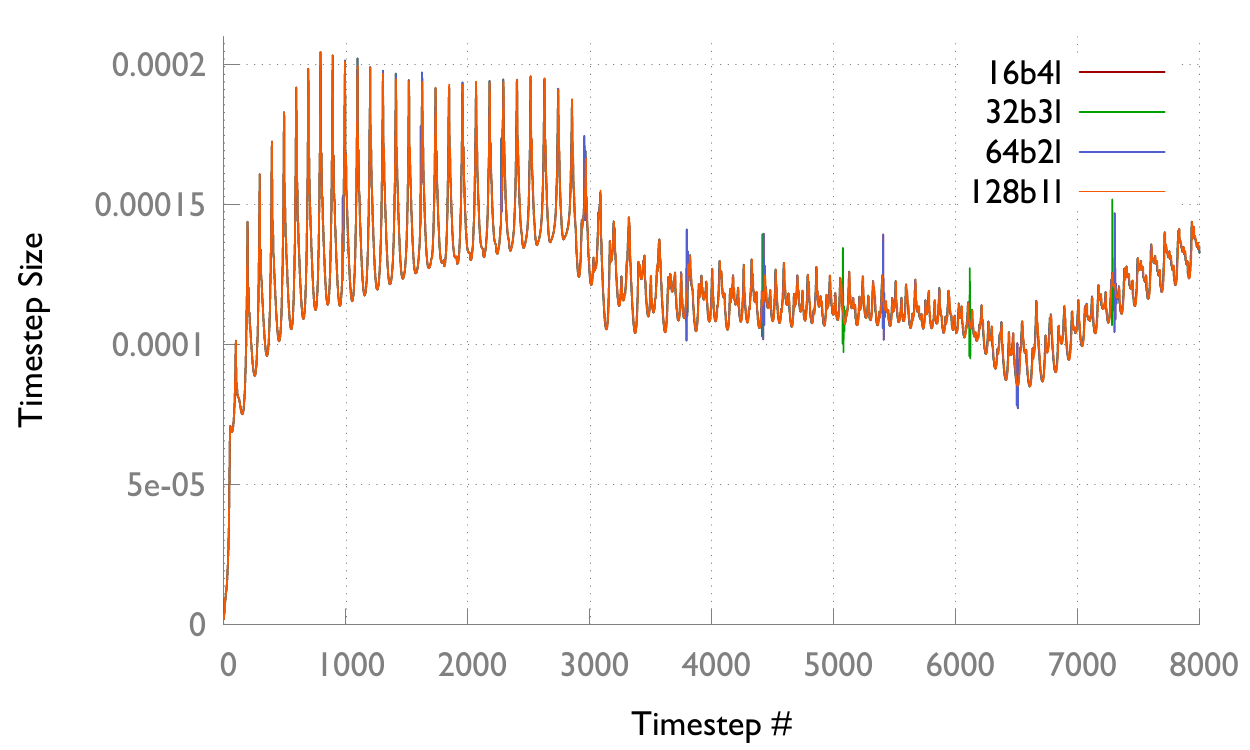} 
\caption{Timestep evolution for different AMR grids with equivalent resolution}
\label{fig:dthistory}
\end{figure}

Table \ref{table:ndt} presents the total number of timesteps taken for each simulation. Little or no variation is seen moving diagonally across the table from left to right. Figure \ref{fig:dthistory} plots the per step time step variation over the course of the simulations on different equivalent resolution AMR grids. Other than at regrid events the figure shows close tracking of the timestep sizes for the various grids. The periodic rise and fall of the timestep size as the simulation progresses is characteristic of truncation error based timestep controllers. From this and the previous tables on linear and nonlinear iteration counts we conclude that for our application the cost per timestep does not rise when moving from uniform grids to AMR grids of equivalent fine resolution.

\begin{figure}[!h]
\begin{center}
\includegraphics[width=0.9\linewidth]{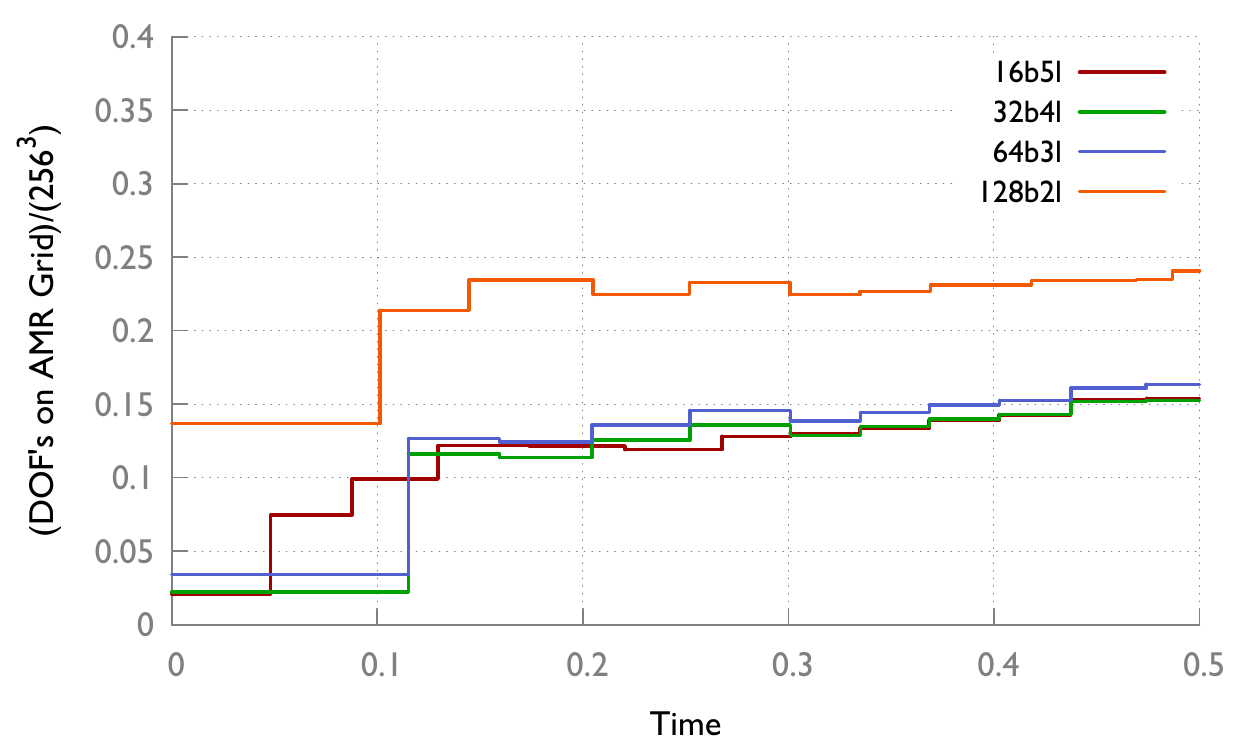}
\end{center}
\caption{Relative degrees of freedom for AMR simulations (relative to uniform $256^3$ simulations)}
\label{case1-rdofs}
\end{figure}

\begin{table}[h]
\centering
\renewcommand{\arraystretch}{1.3}
\begin{tabular}{|c|c|c|c|c|c|}
\hline
   Levels      & 1 & 2 & 3 & 4 & 5 \\
\hline\hline
$ 16^3      $ & -         & 1333 & 3251 & 8204 & 8697  \\ \hline
$ 32^3      $ & 1334 & 3251 & 8204 & 8704 & -          \\ \hline
$ 64^3      $ & 3253 & 8204 & 8704 & -         & -          \\ \hline
$ 128^3    $ & 8207 & 8734 & -        & -          & -          \\ \hline
$ 256^3    $ & 8698         & -         & -        & -          & -          \\ \hline
\end{tabular}
\caption{\label{table:ndt} Total number of timesteps}
\end{table}%
Having shown the similar performance of the AMR calculations to uniform grid calculations with respect to the convergence of the solvers at each timestep and the total number of time steps we turn our attention to the potential performance gains with AMR. Figure \ref{case1-rdofs} plots the relative degrees of freedom (DOFs) required for AMR grid simulations relative to a uniform $256^3$ grid calculation over several regrid events. We see that introducing one level of refinement and decreasing the coarse grid resolution to $128^3$ reduces the number of required degrees of freedom significantly to less than $25\%$ of a uniform grid calculation. As can be seen further reductions in the number of degrees of freedom are obtained as the coarse grid resolution is reduced and more AMR levels are introduced. We note that the reduced degrees of freedom required for an AMR calculation have a significant benefit on the memory requirements for Newton-Krylov methods as the memory required for storing the Krylov subspace vectors is significantly reduced. Combined with strong multilevel preconditioners that reduce the size of the Krylov subspace required at each Newton step we significantly lower the memory requirements for our calculations by using AMR.

\begin{figure}[h]
\begin{center}
\includegraphics{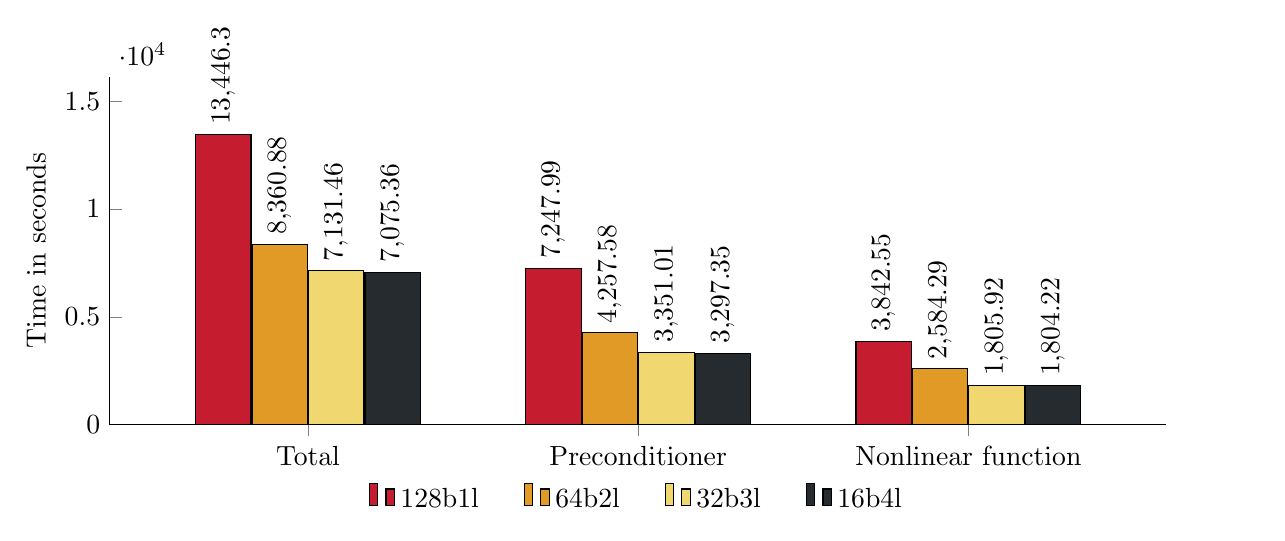}
\caption{Timings for 256b1l equivalent AMR grids on 128 cores}
\label{case1-wallclock}
\end{center}
\end{figure}

The data presented on the number of linear and nonlinear iterations and the total number of timesteps taken indicates that solution efficiency in terms of convergence rates and minimizing the total number of timesteps required to integrate to a given time can be maintained by strong preconditioning and modifications to the timestep control strategy after regridding. Figure \ref{case1-wallclock} presents the wall clock time required to integrate to time $t=1.0$ for the grids shown on 128 cores of an ORNL Cray XK6 machine. We note that while the code has been profiled and some optimizations done further code optimizations are actively being pursued. In particular, preliminary profiling of our application indicates that there is significant room to optimize the communication intensive preconditioning phases of the simulation. Furthermore, based on Figure \ref{case1-rdofs} we can conclude that 128 cores are probably not required to run the grids with refinement levels as  opposed to the uniform $256^3$ grid calculation as on average less than 20\% DOF's are required to perform the AMR calculations. Nevertheless, keeping the number of cores fixed we show that the AMR calculations do indeed significantly reduce the wall clock times required. We note that a significant gain is obtained when introducing one and two levels of refinement  while the improvements are less marked after that point. This is consistent with the data in the previous table where we note that higher levels of refinement did not yield much efficiency. We believe these numbers can potentially be improved by introduced improvements to our refinement tagging strategies and by optimizing parallel performance. We hope to report on that work in future.

\subsection{Accuracy study}

\begin{table}[h]
\centering
\renewcommand{\arraystretch}{1.3}
\begin{tabular}{|c|c|c|c|c|c|}
\hline
$t$ & 0.05 & 0.15 & 0.25 & 0.35 & 0.45 \\ \hline\hline
$\Delta t$ & \multicolumn{5}{|c|}{$L2$ Error: Energy Density} \\ \hline
2.0e-04& 4.53e-05  & 3.93e-05  & 2.87e-05  & 2.47e-05  & 2.36e-05 \\ \hline
1.0e-04& 1.13e-05  & 9.80e-06  & 7.10e-06  & 6.10e-06  & 5.80e-06 \\ \hline
5.0e-05& 2.30e-06  & 2.00e-06  & 1.50e-06  & 1.30e-06  & 1.20e-06 \\ \hline
 & \multicolumn{5}{|c|}{$L2$ Error: Temperature} \\ \hline
2.0e-04& 2.43e-05  & 2.06e-05  & 1.89e-05  & 1.87e-05  & 1.82e-05 \\ \hline
1.0e-04& 6.00e-06  & 5.10e-06  & 4.70e-06  & 4.60e-06  & 4.50e-06 \\ \hline
5.0e-05& 1.20e-06  & 1.00e-06  & 1.00e-06  & 9.00e-07  & 9.00e-07 \\ \hline
\end{tabular}
\caption{\label{table:temporal error} Temporal $L2$ norm errors on a 16b4l grid}
\end{table}%
The previous section focussed on the performance gains
obtained with AMR. In this section we will document the temporal
and spatial accuracy for our simulations on AMR grids.
%
%\begin{figure}[!ht]
%\centering
%\input{plots/tEnergyError.tex}
%\caption{Temporal energy errors on the 128b1l (left)
%and 16b4l (right) grids.}
%\label{fig:tEnergyError}
%\end{figure}
%
%\begin{figure}[!ht]
%\centering
%\input{plots/tTempError.tex}
%\caption{Temporal temperature errors on the 128b1l (left)
%and 16b4l (right) grids.}
%\label{fig:tTempError}
%\end{figure}

\noindent {\it Temporal accuracy: }
In the absence of analytic solutions we first compute reference solutions at different temporal data points in the time interval $[0,0.5]$ on a uniform 128b1l grid with a fixed timestep of $ 2.5 \times 10^{-5}$. 
While maintaining the same equivalent spatial fine grid resolution to keep the contribution from the spatial discretization error the same, simulations
with fixed timesteps of $ \Delta t = 2 \times 10^{-4} $, $ 1 \times 10^{-4} $ and $ 5 \times 10^{-5} $ are run and the solutions at the same data points are compared to the reference solution.
Table \ref{table:temporal error} presents the L2 norm errors for the computed energy density and temperature on 16b4l AMR grids
when compared against the reference solution. The errors decrease by a factor of at least 4 when the timestep is halved showing the second order accuracy of the BDF2 time integration scheme. While not shown, for a given timestep the differences between the solutions on uniform and AMR grids with the same equivalent finest grid resolution are also significantly lower than the spatial and temporal discretization errors.

\begin{table}[h]
\centering
\renewcommand{\arraystretch}{1.3}
\begin{tabular}{|c|c|c|c|c|c|c|}
\hline
$t$ & 0.05 & 0.11 & 0.19 & 0.27 & 0.36 & 0.45 \\ \hline\hline
Grid & \multicolumn{6}{|c|}{$L2$ Error: Energy Density} \\ \hline
16b1l & 1.22e-01  & 3.46e-01  & 4.45e-01  & 4.65e-01  & 4.22e-01  & 3.77e-01 \\ \hline
16b2l & 2.28e-01  & 1.99e-01  & 1.62e-01  & 1.47e-01  & 1.35e-01  & 1.13e-01 \\ \hline
16b3l & 9.27e-02  & 5.64e-02  & 5.09e-02  & 4.92e-02  & 3.20e-02  & 4.13e-02 \\ \hline
16b4l & 1.94e-02  & 1.30e-02  & 1.08e-02  & 1.07e-02  & 1.06e-02  & 9.15e-03 \\ \hline 
%32b1l & 2.28e-01  & 1.99e-01  & 1.62e-01  & 1.47e-01  & 1.35e-01  & 1.13e-01 \\ \hline
%64b1l & 9.27e-02  & 5.64e-02  & 5.09e-02  & 4.92e-02  & 3.20e-02  & 4.13e-02 \\ \hline
128b1l & 1.94e-02  & 1.30e-02  & 1.08e-02  & 1.07e-02  & 1.06e-02  & 9.12e-03 \\ \hline \hline
 & \multicolumn{6}{|c|}{$L2$ Error: Temperature} \\ \hline
16b1l & 5.57e-02  & 1.61e-01  & 2.21e-01  & 2.30e-01  & 2.34e-01  & 2.29e-01 \\ \hline
16b2l & 8.03e-02  & 8.05e-02  & 7.54e-02  & 7.98e-02  & 7.04e-02  & 6.12e-02 \\ \hline
16b3l & 2.55e-02  & 2.28e-02  & 2.01e-02  & 2.17e-02  & 1.69e-02  & 1.99e-02 \\ \hline
16b4l & 5.15e-03  & 4.37e-03  & 3.91e-03  & 3.81e-03  & 4.14e-03  & 4.20e-03 \\ \hline 
%32b1l & 8.03e-02  & 8.05e-02  & 7.54e-02  & 7.98e-02  & 7.04e-02  & 6.12e-02 \\ \hline
%64b1l & 2.55e-02  & 2.28e-02  & 2.01e-02  & 2.17e-02  & 1.69e-02  & 1.99e-02 \\ \hline
128b1l & 5.15e-03  & 4.37e-03  & 3.91e-03  & 3.81e-03  & 4.14e-03  & 4.19e-03 \\ \hline
\end{tabular}
\caption{\label{table:spatialerror} Spatial $L2$ norm errors}
\end{table}%
%
%\begin{figure}[!ht]
%\centering
%\input{plots/sEnergyError.tex}
%\caption{Spatial energy errors on the uniform (left)
%and AMR (right) grids.}
%\label{fig:sEnergyError}
%\end{figure}
%
%\begin{figure}[!ht]
%\centering
%\input{plots/sTempError.tex}
%\caption{Spatial temperature errors on the 128b1l (left)
%and 16b4l (right) grids.}
%\label{fig:sTempError}
%\end{figure}
%
\noindent {\it Spatial accuracy:} 
For the spatial accuracy studies, as before, since an exact solution is not available,
a simulation is done on a uniform $ 256^3 $ grid and the time steps are recorded as well
as the solution at different points in time. Solutions obtained using the same timesteps on AMR grids (to account for the effect of temporal discretization error) are then compared to the
reference solutions on the uniform reference grid. Table \ref{table:spatialerror}
presents the L2 norm of the spatial discretization error
for energy density and temperature respectively . Note that these runs were single material runs with
the atomic number $ z = 1 $ in the whole domain.
Second order accuracy is obtained though some loss of accuracy is observed in the initial stages of the calculation.
For reference the errors on a uniform 128b1l grid are also shown in the last row for the energy density and temperature respectively.
\begin{table}[h]
\centering
\renewcommand{\arraystretch}{1.3}
\begin{tabular}{|c|c|c|c|c|c|c|}
\hline
$t$ & 0.04 & 0.10 & 0.19 & 0.25 & 0.35 & 0.46 \\ \hline\hline
Grid & \multicolumn{6}{|c|}{$L2$ Error: Energy Density} \\ \hline
16b1l & 6.01e-02  & 3.32e-01  & 4.55e-01  & 4.99e-01  & 4.49e-01  & 3.77e-01 \\ \hline
16b2l & 2.21e-01  & 2.12e-01  & 1.63e-01  & 1.55e-01  & 1.28e-01  & 1.15e-01 \\ \hline
16b3l & 9.92e-02  & 6.82e-02  & 4.67e-02  & 4.12e-02  & 3.84e-02  & 4.28e-02 \\ \hline
16b4l & 2.38e-02  & 1.39e-02  & 1.17e-02  & 1.11e-02  & 1.44e-02  & 2.29e-02 \\ \hline \hline
%32b1l & 2.21e-01  & 2.12e-01  & 1.63e-01  & 1.55e-01  & 1.28e-01  & 1.15e-01 \\ \hline
%64b1l & 9.92e-02  & 6.82e-02  & 4.67e-02  & 4.12e-02  & 3.84e-02  & 4.28e-02 \\ \hline
128b1l & 2.38e-02  & 1.39e-02  & 1.17e-02  & 1.11e-02  & 1.44e-02  & 2.29e-02 \\ \hline \hline
 & \multicolumn{6}{|c|}{$L2$ Error: Temperature} \\ \hline
16b1l & 3.15e-02  & 1.52e-01  & 2.21e-01  & 2.32e-01  & 2.23e-01  & 2.06e-01 \\ \hline
16b2l & 7.78e-02  & 8.97e-02  & 7.55e-02  & 7.17e-02  & 7.08e-02  & 7.86e-02 \\ \hline
16b3l & 2.76e-02  & 2.26e-02  & 2.09e-02  & 1.90e-02  & 2.65e-02  & 3.39e-02 \\ \hline
16b4l & 5.53e-03  & 4.43e-03  & 5.45e-03  & 7.63e-03  & 1.11e-02  & 1.64e-02 \\ \hline \hline
%32b1l & 7.78e-02  & 8.97e-02  & 7.55e-02  & 7.17e-02  & 7.08e-02  & 7.86e-02 \\ \hline
%64b1l & 2.76e-02  & 2.26e-02  & 2.09e-02  & 1.90e-02  & 2.65e-02  & 3.39e-02 \\ \hline
128b1l & 5.53e-03  & 4.43e-03  & 5.45e-03  & 7.63e-03  & 1.11e-02  & 1.64e-02 \\ \hline
\end{tabular}
\caption{\label{table:spatialerror2m} Spatial $L2$ norm errors}
\end{table}
Table \ref{table:spatialerror2m} presents accuracy studies for simulations on the physical domain as shown in Figure \ref{fig:material} on both AMR and a uniform grid (128b1l).
The energy and temperature errors show second order accuracy before time 0.25. However, after time 0.25 the front of the Marshak wave hits a discontinuity across material interfaces and the accuracy drops. While not shown this drop in accuracy is seen on uniform grids also in our tests and appears to be related to the spatial discretization of the nonlinear diffusion coefficients across discontinuous interfaces.

\section{Conclusions and Future Directions}

In this paper we described research on an efficient solution methodology for solving 3D non-equilibrium radiation diffusion problems. Implicit time integration enabled time stepping at the dynamical timescale of the problem. Control theory based step size control minimized the overall required number of steps while allowing us to use methods that have computational stability. Inexact Newton methods as implemented in the JFNK solver minimized the work required at the outer Newton iteration for each time step with GMRES providing a robust solver for non-symmetric definite systems at each Newton iteration. The multilevel preconditioner components were critical in provided level independent convergence of the linear solver. 3D AMR minimized the computational and memory requirements at each step of the calculation. The individual techniques described are not new, but their combination and application to solving problems in radiation transport is new to our knowledge. 

Care had to be taken in selecting and combining the individual components so that the overall simulation was accurate and efficient. Our experience in developing efficient simulation methods for this application and more broadly time dependent nonlinear multiphysics systems is that it requires not only focussing on how the individual simulation components can be efficient and/or accurate but also on understanding how the interplay between the components can enhance or degrade the efficiency and accuracy of the overall simulation methodology. This is true for the interplay between the time step control algorithm and regridding for AMR, time step control and the accuracy of the nonlinear solvers, and controlling the efficiency of the Krylov solvers through level independent preconditioner components to mention a few.

In future work we hope to report on the performance of this application on hybrid multicore-GPU petascale platforms such as the Titan supercomputer at Oak Ridge National Laboratory. The non-equilibrium radiation diffusion application is an excellent testbed for studying the performance of nonlinear multi-physics application solution components on AMR grids at the petascale. Enabling greater asynchrony in our solver components by using algorithms such as AFACx \cite{lee2003asynchronous,lee2004asynchronous}, load balancing for AMR applications, a posteriori error estimation for finite volume AMR applications, e.g. \cite{estep2009posteriori}, and multiphysics smoother components tuned for hybrid architectures are some of the areas we hope to make progress on. 
\section{Acknowledgements}
The authors are grateful to the National Center for Computational Science at ORNL for access to their internal high performance computing resources and to James Schwarzmeier from 
CRAY Inc. for providing access to a CRAY XE6 system for extensive testing of the non-equilibrium radiation diffusion application presented in this paper. Zhen Wang would like to acknowledge support from the Mathematical, Information, and Computational Sciences Division, Office of Advanced Scientific Computing Research, U.S. Department of Energy, under Contract No. DE-AC05-00OR22725 with UT-Battelle, LLC. 
Mark Berrill would like to acknowledge support from the Eugene P. Wigner Fellowship at Oak Ridge National Laboratory, managed by UT-Battelle, LLC, for the U.S. Department of Energy 
under Contract DE-AC05-00OR22725.

\appendix

\section{Figures}

\begin{figure}[ht]
\centering
\includegraphics{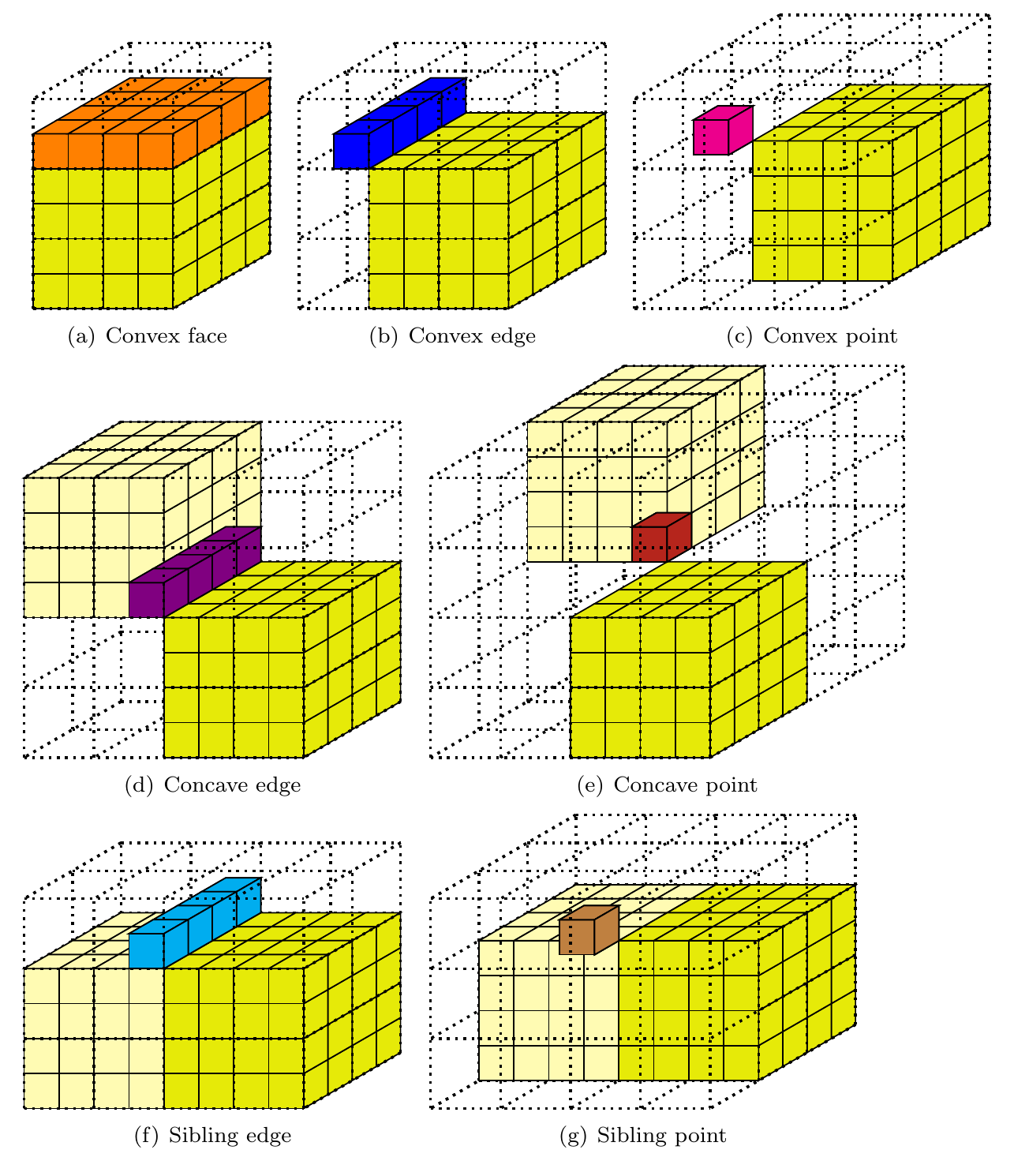}
\caption{Different types of coarse fine boundary fragments
of the darker patch.}
\label{CoarseFineBoundaryFragments}
\end{figure}

\begin{figure}[ht]
\centering
\includegraphics{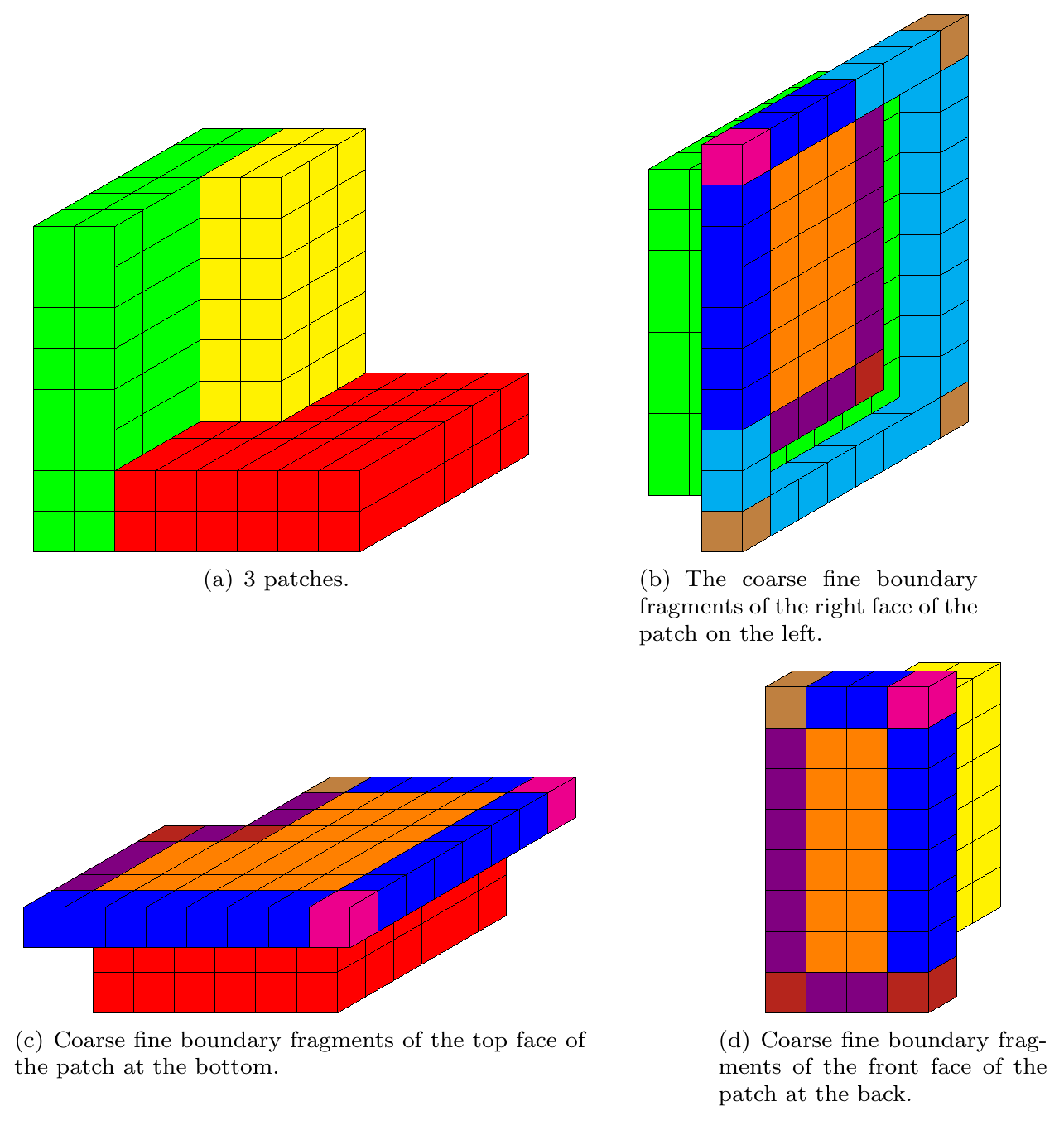}
\caption{A complex example showing different coarse fine boundary fragments.}
\label{CoarseFineBoundaryExample}
\end{figure}

\bibliographystyle{elsarticle-num}
\bibliography{nrdf}

\end{document}